\documentclass{article}
\usepackage{amsfonts}

\usepackage{amsmath}

\newtheorem{theorem}{Theorem}[section]

\newtheorem{corollary}[theorem]{Corollary}

\newtheorem{definition}[theorem]{Definition}

\newtheorem{lemma}[theorem]{Lemma}

\newtheorem{proposition}[theorem]{Proposition}

\newenvironment{proof}[1][Proof]{\textbf{#1.} }{\ \rule{0.5em}{0.5em}}

\begin{document}

\title{\textbf{Rough Path Analysis Via Fractional Calculus}}
\author{Yaozhong Hu\thanks{%
Y. Hu is supported in part by the National Science Foundation under Grant
No. DMS0204613 and DMS0504783 } \ \ and\ 
  David Nualart\thanks{%
The work of D. Nualart is partially supported by the MCyT Grant BFM2000-0598} \\
Department of Mathematics\thinspace ,\ University of Kansas\\
405 Snow Hall\thinspace ,\ Lawrence, Kansas 66045-2142\\
}
\date{}
\maketitle

\begin{abstract}
Using fractional calculus we define integrals of the form $%
\int_{a}^{b}f(x_{t})dy_{t}$, where $x$ and $y$ are vector-valued H\"{o}lder continuous
functions of order $\displaystyle \beta \in (\frac13 ,\  \frac12)$ and $f$ is a 
continuously differentiable function
such that $f'$ is $\lambda$-H\"oldr continuous for some $\lambda>\frac1\beta-2$. 
Under some further smooth conditions on $f$ the integral is a
continuous functional of $x$, $y$,  and the tensor product $x\otimes y$ with
respect to the H\"{o}lder norms. We derive some estimates for these
integrals and we solve differential equations driven by the function $y$. We
discuss some applications to stochastic integrals and stochastic
differential equations.
\end{abstract}

\section{Introduction}

The theory of rough path analysis has been developed from the seminal paper
by Lyons \ \cite{lyons}. The purpose of this theory is to analyze dynamical
systems $dx_{t}=f(x_{t})dy_{t}$, where the control function $y$ is not
differentiable. If the rough control $y$ has finite $p$-variation on bounded
intervals, where $p\geq 2$, then the dynamical system is a continuous
function, in the $p$-variation norm, of $y$ and the associated
multiplicative functionals $\overbrace{y\otimes {\cdots }\otimes y}^k $, with $%
k=2,\ldots ,[p]$. In the case $1\leq p<2$, the dynamical system can be
formulated using Riemann-Stieltjes integrals and applying the results of
Young \cite{Yo}. In this case, $x_{t}$ is a continuous function of $y$ in
the $p$-variation norm (see Lyons \cite{Ly2}).

Suppose that $f$ and $g$ are H\"{o}lder continuous functions on the interval 
$[a,b]$, of order $\lambda $ and $\mu $, respectively, with $\lambda +\mu >1$%
. Then, the Riemann-Stieltjes integral $\int_{a}^{b}fdg$ can be expressed as
a Lebesgue integral using fractional derivatives (see Z\"{a}hle \cite{Za98}
and Proposition 2.1 below). This fact has been exploited by Nualart and R%
\u{a}\c{s}canu in \cite{NR} to analyze dynamical systems driven by a control
function $y$ which is H\"{o}lder continuous of order $\beta >\frac{1}{2}$.
In this case further results are obtained in \cite{HN} along 
the line  of the present  paper. 

The purpose of this paper is to analyze dynamical systems $%
dx_{t}=f(x_{t})dy_{t}$, \ where the control function $y$ is H\"{o}lder
continuous of order $\beta \in (\frac{1}{3},\frac{1}{2})$, using the
techniques of the classical fractional calculus, and  following an approach
inspired by the work of Nualart and R\u{a}\c{s}canu \cite{NR} in the case $%
\beta >\frac{1}{2}$. In order to achieve this objective, we first provide in
Section 3 an explicit formula for integrals of the form $%
\int_{a}^{b}f(x_{t})dy_{t}$, where $x$ and $y$ are H\"{o}lder continuous of
order $\beta \in (\frac{1}{3},\frac{1}{2})$. This formula, given in Theorem
3.1, is based on the fractional integration by parts formula, and it
involves the functions $x$, $y$, and the quadratic multiplicative functional 
$x\otimes y\,$. Notice that this explicit formula does not depend on any
approximation scheme. \ As a consequence, we derive estimates in the H\"{o}%
lder norm for the indefinite integral.

Section 4 is devoted to establish the existence and uniqueness of a solution
for the dynamical system $dx_{t}=f(x_{t})dy_{t}$. \ The main ingredient in
the proof of these results is to transform this equation into a system of
integral equations for $x$ and $x\otimes y$ that can be solved by a standard
application of a fixed point argument. We  show how  the solution
depends continuously on the H\"{o}lder norm of $y$ and $y\otimes y$.
We also prove some stability results for the differential equations
which are interesting, new and may be difficult to obtain by
other approaches. 
Remark that to derive our results we do not make use of the theory of rough
paths, and we obtain explicit formulas that do not depend on any
approximation argument.

These results can be applied to implement a path-wise approach to define
stochastic integrals and to solve stochastic differential equations driven
by a multidimensional Brownian motion.  As an application of the 
deterministic results obtained for dynamical systems we derive a 
sharp rate of almost sure convergence of the Wong-Zakai approximation for 
multidimensional diffusion processes.   We couldn't find
this kind of estimates elsewhere. 
Similar results hold in the case of a fractional Brownian motion 
with Hurst parameter $H\in (\frac{1}{3},\frac{1}{2%
})$.  The approximation
of the solutions of stochastic differential equations driven by a fractional
Brownian motion with Hurst parameter $H\in (\frac{1}{3},\frac{1}{2})$ is
more involved and it will be treated in a forthcoming paper.

\section{Fractional Integrals and Derivatives}

Let $a,b\in \mathbb{R}$ with $a<b.$ Denote by $L^{p}\left( a,b\right) $, $%
p\geq 1$, the space of Lebesgue measurable functions $f:\left[ a,b\right]
\rightarrow \mathbb{R}$ for which $\left\Vert f\right\Vert _{L^{p}\left(
a,b\right) }<\infty $, where 
\begin{equation*}
\left\Vert f\right\Vert _{L^{p}\left( a,b\right) }=\left\{ 
\begin{array}{lll}
\left( \int_{a}^{b}\left\vert f\left( t\right) \right\vert ^{p}dt\right)
^{1/p}, &  & \hbox{if \ }1\leq p<\infty \\ 
ess\sup \left\{ \left\vert f\left( t\right) \right\vert :t\in \left[ a,b%
\right] \right\} , &  & \hbox{if \ }p=\infty .%
\end{array}%
\right.
\end{equation*}

Let $f\in L^{1}\left( a,b\right) $ and $\alpha >0.$ The left-sided and
right-sided fractional Riemann-Liouville integrals of $f$ of order $\alpha $
are defined for almost all $x\in \left( a,b\right) $ by 
\begin{equation*}
I_{a+}^{\alpha }f\left( t\right) =\frac{1}{\Gamma \left( \alpha \right) }%
\int_{a}^{t}\left( t-s\right) ^{\alpha -1}f\left( s\right) ds
\end{equation*}%
and 
\begin{equation*}
I_{b-}^{\alpha }f\left( t\right) =\frac{\left( -1\right) ^{-\alpha }}{\Gamma
\left( \alpha \right) }\int_{t}^{b}\left( s-t\right) ^{\alpha -1}f\left(
s\right) ds,
\end{equation*}%
respectively, where $\left( -1\right) ^{-\alpha }=e^{-i\pi \alpha }$ and $%
\Gamma \left( \alpha \right) =\int_{0}^{\infty }r^{\alpha -1}e^{-r}dr$ is
the Euler gamma function. Let $I_{a+}^{\alpha }(L^{p})$ (resp. $%
I_{b-}^{\alpha }(L^{p})$) be the image of $L^{p}(a,b)$ by the operator $%
I_{a+}^{\alpha }$ (resp. $I_{b-}^{\alpha }$). If $f\in I_{a+}^{\alpha
}\left( L^{p}\right) \ $ (resp. $f\in I_{b-}^{\alpha }\left( L^{p}\right) $)
and $0<\alpha <1$ then the Weyl derivatives are defined as 
\begin{equation}
D_{a+}^{\alpha }f\left( t\right) =\frac{1}{\Gamma \left( 1-\alpha \right) }%
\left( \frac{f\left( t\right) }{\left( t-a\right) ^{\alpha }}+\alpha
\int_{a}^{t}\frac{f\left( t\right) -f\left( s\right) }{\left( t-s\right)
^{\alpha +1}}ds\right)  \label{1.1}
\end{equation}%
and 
\begin{equation}
D_{b-}^{\alpha }f\left( t\right) =\frac{\left( -1\right) ^{\alpha }}{\Gamma
\left( 1-\alpha \right) }\left( \frac{f\left( t\right) }{\left( b-t\right)
^{\alpha }}+\alpha \int_{t}^{b}\frac{f\left( t\right) -f\left( s\right) }{%
\left( s-t\right) ^{\alpha +1}}ds\right)  \label{1.2}
\end{equation}%
where $a\leq t\leq b$ (the convergence of the integrals at the singularity $%
s=t$ holds point-wise for almost all $t\in \left( a,b\right) $ if $p=1$ and
moreover in $L^{p}$-sense if $1<p<\infty $).

For any $\lambda \in (0,1)$, we denote by $C^{\lambda }(a,b)$ the space of $%
\lambda $-H\"{o}lder continuous functions on the interval $[a,b]$. Recall
from \cite{SaKiMa93} that we have:

\begin{itemize}
\item If $\alpha <\frac{1}{p}$ and $q=\frac{p}{1-\alpha p}$ then 
\begin{equation*}
I_{a+}^{\alpha }\left( L^{p}\right) =I_{b-}^{\alpha }\left( L^{p}\right)
\subset L^{q}\left( a,b\right) .
\end{equation*}

\item If $\alpha >\frac{1}{p}$ then%
\begin{equation*}
I_{a+}^{\alpha }\left( L^{p}\right) \,\cup \,I_{b-}^{\alpha }\left(
L^{p}\right) \subset C^{\alpha -\frac{1}{p}}\left( a,b\right) .
\end{equation*}
\end{itemize}

The following inversion formulas hold: 
\begin{eqnarray}
I_{a+}^{\alpha }\left( D_{a+}^{\alpha }f\right) &=&f,\quad \quad \;\forall
f\in I_{a+}^{\alpha }\left( L^{p}\right)  \label{1.4} \\
I_{a-}^{\alpha }\left( D_{a-}^{\alpha }f\right) &=&f,\quad \quad \;\forall
f\in I_{a-}^{\alpha }\left( L^{p}\right)  \label{1.3}
\end{eqnarray}%
and 
\begin{equation}
D_{a+}^{\alpha }\left( I_{a+}^{\alpha }f\right) =f,\quad D_{a-}^{\alpha
}\left( I_{a-}^{\alpha }f\right) =f,\quad \;\forall f\in L^{1}\left(
a,b\right) \,.  \label{1.5}
\end{equation}%
On the other hand, for any $f,g\in L^{1}(a,b)$ we have 
\begin{equation}
\int_{a}^{b}I_{a+}^{\alpha }f(t)g(t)dt=(-1)^{\alpha
}\int_{a}^{b}f(t)I_{b-}^{\alpha }g(t)dt\,,  \label{1.6}
\end{equation}%
and for $f\ \in I_{a+}^{\alpha }\left( L^{p}\right) $ and $g\in
I_{a-}^{\alpha }\left( L^{p}\right) $ we have%
\begin{equation}
\int_{a}^{b}D_{a+}^{\alpha }f(t)g(t)dt=(-1)^{-\alpha
}\int_{a}^{b}f(t)D_{b-}^{\alpha }g(t)dt.  \label{1.7}
\end{equation}

Suppose that $f\in C^{\lambda }(a,b)$ and $g\in C^{\mu }(a,b)$ with $\lambda
+\mu >1$. Then, from the classical paper by Young \cite{Yo}, the
Riemann-Stieltjes integral $\int_{a}^{b}fdg$ exists. The following
proposition can be regarded as a fractional integration by parts formula,
and provides an explicit expression for the integral $\int_{a}^{b}fdg$ in
terms of fractional derivatives (see \cite{Za98}).

\begin{proposition}
\label{p.2.1} Suppose that $f\in C^{\lambda }(a,b)$ and $g\in C^{\mu }(a,b)$
with $\lambda +\mu >1$. Let ${\lambda }>\alpha $ and $\mu >1-\alpha $. Then
the Riemann Stieltjes integral $\int_{a}^{b}fdg$ exists and it can be
expressed as%
\begin{equation}
\int_{a}^{b}fdg=(-1)^{\alpha }\int_{a}^{b}D_{a+}^{\alpha }f\left( t\right)
D_{b-}^{1-\alpha }g_{b-}\left( t\right) dt,  \label{1.8}
\end{equation}%
where $g_{b-}\left( t\right) =g\left( t\right) -g\left( b\right) $.
\end{proposition}

\bigskip

We will make use of the following two-variable fractional integration by
parts formula, whose proof is given in the Appendix.

\begin{lemma}
\label{l.2.2} Let $\varphi (\xi ,\eta )$ and $\psi (\xi ,\eta )$ be two
functions of class $C^{2}$ defined on $a\leq \xi \leq \eta \leq b$. Suppose $%
\psi (\xi ,\eta )$ vanishes on the diagonal. The following fractional
integration by parts formula holds for any $0<\alpha <1$.%
\begin{equation}
\int_{a}^{b}d\xi \int_{\xi }^{b}\varphi (\xi ,\eta )\frac{\partial ^{2}\psi 
}{\partial \xi \partial \eta }(\xi ,\eta )d\xi d\eta =-\int_{a}^{b}d\eta
\int_{a}^{\eta }D_{a+}^{\alpha ,\xi }D_{b-}^{\alpha ,\eta }\varphi
_{a+,b-}(\xi ,\eta )\Gamma ^{\alpha }\psi (\xi ,\eta )d\xi ,  \label{1.9}
\end{equation}%
where $D_{a+}^{\alpha ,\xi }$ denotes the fractional derivative on variable $%
\xi $, $D_{b-}^{\alpha ,\eta }$ denotes the fractional derivative on the
variable $\eta $, and the operator $\Gamma ^{\alpha }$ is defined by%
\begin{equation}
\Gamma ^{\alpha }\psi (\xi ,\eta )=D_{\eta -}^{1-\alpha ,\xi }D_{\xi
+}^{1-\alpha ,\eta }\psi (\xi ,\eta ).  \label{1.10}
\end{equation}
\end{lemma}

 \setcounter{equation}{0}
\section{Integration of Rough Functions}

Fix $\frac{1}{3}<\beta <\frac{1}{2}$. Suppose that $x\ :[0,T]\rightarrow 
\mathbb{R}^{m}$ and $y\ :[0,T]\rightarrow \mathbb{R}^{d}$ are $\beta $-H\"{o}%
lder continuous functions. Following \cite{lyons} we assume that $x\otimes y$
is well-defined and $\ $it is a continuous function defined on $\Delta
:=\{(s,t):0\leq s\leq t\leq T\}$ with values on $\mathbb{R}^{m}\otimes 
\mathbb{R}^{d}$ verifying the following properties:

\begin{description}
\item[i)] For all $s\leq u\leq t$ we have (\textit{multiplicative property})%
\begin{equation}
\left( x\otimes y\right) _{s,u}+\left( x\otimes y\right)
_{u,t}+(x_{u}-x_{s})\otimes (y_{t}-y_{u})=\left( x\otimes y\right) _{s,t}.
\label{2.1}
\end{equation}

\item[ii)] For all $(s,t)\in \Delta $%
\begin{equation}
\left| \left( x\otimes y\right) _{s,t}\right| \leq k|t-s|^{2\beta }.
\label{2.2}
\end{equation}
\end{description}

That is, $(x,y,x\otimes y)$ constitutes a multiplicative functional in the
sense of the rough paths analysis theory. We will say that $(x,y,x\otimes y)$
is a $\beta $-\textit{H\"{o}lder continuous multiplicative functional on }$%
\mathbb{R}^{m}\otimes \mathbb{R}^{d}$.

If $x$ and $y$ are smooth functions, then 
\begin{equation}
\left( x\otimes y\right) _{s,t}^{i,j}=\int_{s<\xi <\eta <t}dx_{\xi
}^{i}dy_{\eta }^{j}  \label{2.3}
\end{equation}%
clearly defines a $\beta $-H\"{o}lder continuous multiplicative functional.

Let $f:\mathbb{R}^{m}\mathbb{\rightarrow R}^{d}$ be a continuously
differentiable function such that $f^{\prime }$ is $\lambda $-H\"{o}lder
continuous, where $\lambda >\frac{1}{\beta }-2$. Our aim is to define the
integral%
\begin{equation}
\int_{a}^{b}f(x_{r})dy_{r}=\sum_{i=1}^{d}\int_{a}^{b}f_{i}(x_{r})dy_{r}^{i}
\label{2.4}
\end{equation}%
using fractional calculus.

Fix a number $\alpha $ such that $1-\beta <\alpha <2\beta $ and $\alpha <%
\frac{\lambda \beta +1}{2}$. This is possible because $3\beta >1$ and $\frac{%
\lambda \beta +1}{2}>1-\beta $.

Notice first that the fractional integration by parts formula (\ref{1.8})
cannot be used to define the integral (\ref{2.4}) because the fractional
derivative $D_{a+}^{\alpha }f\left( x\right) $ is not well-defined under our
hypotheses. For this reason we introduce the following \textit{compensated
fractional derivative}:

\begin{eqnarray}
\lefteqn{\widehat{D}_{a+}^{\alpha }f\left( x\right) (r)=\frac{1}{\Gamma \left(
1-\alpha \right) }\Bigg( \frac{f\left( x_{r}\right) }{(r-a)^{\alpha }}} \notag  \\%
&&+\alpha \int_{a}^{r}\frac{f\left( x_{r}\right) -f\left( x_{\theta }\right)
-f^{\prime }(x_{\theta })(x_{r}-x_{\theta })}{\left( r-\theta \right)
^{\alpha +1}}d\theta \Bigg) .  \label{2.5}
\end{eqnarray}%
This derivative is well-defined under our hypotheses because%
\begin{equation*}
\frac{\left| f\left( x_{r}\right) -f\left( x_{\theta }\right) -f^{\prime
}(x_{\theta })(x_{r}-x_{\theta })\right| }{\left( r-\theta \right) ^{\alpha
+1}}\leq K|r-\theta |^{(1+\lambda )\beta -\alpha -1},
\end{equation*}%
where $K=\left\| f^{\prime }\right\| _{\lambda }\left\| x\right\| _{\beta
}^{1+\lambda }$ and $(1+\lambda )\beta -\alpha >0$ since $\alpha <\frac{%
\lambda \beta +1}{2}<(1+\lambda )\beta $.

For $\theta <\xi <\eta $ introduce the kernel 
\begin{eqnarray}
G(\theta ,\xi ,\eta ):= &&\frac{\newline
1}{\alpha \Gamma \left( \alpha \right) \Gamma (2\alpha -1)}\left( \xi
-\theta \right) ^{\alpha -1}\left( \eta -\xi \right) ^{\alpha -1}  \notag \\
&&\times \int_{0}^{1}q^{2\alpha -2}(1-q)^{-\alpha }\left( 1+(1-q)\frac{\xi
-\theta }{\eta -\xi }\right) ^{-1}dq.  \label{2.6}
\end{eqnarray}%
Define for $\varepsilon <\alpha +\beta -1$ and $\theta <\xi <\eta <b$ 
\begin{equation}
K_{\theta ,b}(\xi ,\eta )=-D_{\theta +}^{1,\alpha -{\varepsilon }%
}D_{b-}^{2,\alpha -{\varepsilon }}\left[ G_{b-}(\theta ,\xi ,\eta )\right]
\,.  \label{2.7}
\end{equation}%
In \ Lemma \ref{l.7.2} we will show that \ this kernel satisfies%
\begin{equation*}
\sup_{0\le s<t\le T}\int_{s< \xi < \eta < t}\left| K_{s,t}(\xi ,\eta
)\right| d\xi d\eta <\infty .
\end{equation*}

Finally, we denote 
\begin{equation}
{\Lambda }_{a}^{b}(x\otimes y):=\int_{a}^{b}\int_{a}^{\eta }K_{a,b}(\xi
,\eta )\Gamma ^{\alpha -{\varepsilon }}\left( x\otimes y\right) _{\xi ,\eta
}d\xi d\eta .  \label{2.8}
\end{equation}

We are ready now to define the integral $\int_{a}^{b}f(x_{r})dy_{r}$.

\begin{definition}
\label{def1} Let $(x,y,x\otimes y)$ be a $\beta $-\textit{H\"{o}lder
continuous multiplicative functional }on\textit{\ }$\mathbb{R}^{m}\otimes 
\mathbb{R}^{d}$. Let $f:\mathbb{R}^{m}\mathbb{\rightarrow R}^{d}$ be a
continuously differentiable function such that $f^{\prime }$ is $\lambda $-H%
\"{o}lder continuous, where $\lambda >\frac{1}{\beta }-2$. Fix $\alpha >0$
and $\varepsilon >0$ such that $1-\beta <\alpha <2\beta $, $\alpha <\frac{%
\lambda \beta +1}{2}$ and $\ $ $\varepsilon <\alpha +\beta -1$. Then, for
any $0\leq a<b\leq T$ we define%
\begin{eqnarray}
\int_{a}^{b}f(x_{r})dy_{r} &=&(-1)^{\alpha }\ \sum_{i=1}^{d}\int_{a}^{b}%
\widehat{D}_{a+}^{\alpha }f_{i}\left( x\right) (r)D_{b-}^{1-\alpha
}y_{b-}^{i}(r)dr  \notag \\
&&+\sum_{i=1}^{m}\sum_{j=1}^{d}\int_{a}^{b}D_{a+}^{2\alpha -1}\partial
_{i}f_{j}\left( x\right) (r){\Lambda }_{r}^{b}(x^{i}\otimes y^{j})dr\,.
\label{2.9}
\end{eqnarray}
\end{definition}

Notice that if $y$ is $\beta $-H\"{o}lder continuous, the fractional
derivative $D_{b-}^{1-\alpha }y_{b-}(r)$ is well-defined because 
\begin{equation*}
\frac{\left\vert y_{\sigma }-y_{r}\right\vert \ }{\left( \sigma -r\right)
^{2-\alpha }}\leq \left\Vert y\right\Vert _{\beta }|\sigma -r|^{\beta
+\alpha -2}
\end{equation*}%
and $\beta +\alpha -2>-1$.The following theorem asserts that this definition
is coherent with the classical notion of integral and will allow us to
deduce estimates in the H\"{o}lder norm.

\begin{theorem}
\label{t.3.1} Suppose $y:[0,T]\rightarrow \mathbb{R}^{d}$ is a continuously differentiable
function. Let $x\ :[0,T]\rightarrow \mathbb{R}^{m}$ be a $\beta $-H\"{o}lder
continuous function and let $x\otimes y$ be defined by $\left( x\otimes
y\right) _{s,t}^{i,j}=\int_{s}^{t}\left( x_{\xi }^{i}-x_{s}^{i}\right)
\left( y^{j}\right) _{\xi }^{\prime }d\xi $. Assume that $f$ satisfies the
assumptions of Definition \ref{def1}. Then, the integral $%
\int_{a}^{b}f(x_{r})dy_{r}$ introduced in (\ref{2.9}) coincides with $%
\sum_{i=1}^d\int_{a}^{b}f_i(x_{r})(y^i)_{r}^{\prime }dr$.
\end{theorem}

\begin{proof}
To simplify the proof we take $m=d=1$. From (\ref{1.8}) and (\ref{2.5}) we
get%
\begin{eqnarray*}
\int_{a}^{b}f(x_{r})y_{r}^{\prime }dr &=&(-1)^{\alpha
}\int_{a}^{b}D_{a+}^{\alpha }f\left( x\right) (r)D_{b-}^{1-\alpha
}y_{b-}(r)dr \\
&=&(-1)^{\alpha }\int_{a}^{b}\widehat{D}_{a+}^{\alpha }f\left( x\right)
(r)D_{b-}^{1-\alpha }y_{b-}(r)dr+A_{2},
\end{eqnarray*}%
where%
\begin{eqnarray}
A_{2} &=&\frac{\alpha (-1)^{\alpha }}{\Gamma \left( 1-\alpha \right) }%
\int_{a}^{b}\int_{a}^{r}\frac{f^{\prime }(x_{\theta })(x_{r}-x_{\theta })}{%
\left( r-\theta \right) ^{\alpha +1}}D_{b-}^{1-\alpha }y_{b-}(r)drd\theta 
\notag \\
&=&\frac{\alpha (-1)^{\alpha }}{\Gamma \left( 1-\alpha \right) }%
\int_{a}^{b}f^{\prime }(x_{\theta })\left( \int_{\theta }^{b}\frac{%
x_{r}-x_{\theta }}{\left( r-\theta \right) ^{\alpha +1}}D_{b-}^{1-\alpha
}y_{b-}(r)dr\right) d\theta .  \label{2.10}
\end{eqnarray}%
So, it suffices to show that%
\begin{equation}
A_{2}=\int_{a}^{b}D_{a+}^{2\alpha -1}f^{\prime }\left( x\right) (r){\Lambda }%
_{r}^{b}(x\otimes y)dr.  \label{2.11}
\end{equation}
Formula (\ref{2.11}) should be first proved for $x$ of class $C^1$ and then extended to a
general $\beta $-H\"{o}lder continuous function. Applying\ (\ref{1.4}), (\ref%
{1.6}), and \ (\ref{1.5}) we obtain%
\begin{eqnarray*}
A_{2} &=&\frac{\alpha (-1)^{\alpha }}{{\Gamma }(1-\alpha )}\int_{a}^{b}\
D_{b-}^{1-\alpha }y_{b-}(r)\int_{a}^{r}\frac{f^{\prime }(x_{\theta
})(x_{r}-x_{\theta })}{(r-\theta )^{\alpha +1}}d\theta dr \\
&=&\frac{\alpha (-1)^{3\alpha -1}}{{\Gamma }(1-\alpha )}\int_{a}^{b}\
D_{b-}^{1-\alpha }y_{b-}(r)\int_{a}^{r}D_{a+}^{2\alpha -1}f^{\prime
}(x)(\theta )I_{r-}^{2\alpha -1,\theta }\left( \frac{x_{r}-x_{\theta }}{%
(r-\theta )^{\alpha +1}}\right) d\theta dr \\
&=&\frac{\alpha (-1)^{3\alpha -1}}{{\Gamma }(1-\alpha )}%
\int_{a}^{b}D_{a+}^{2\alpha -1}f^{\prime }(x)(\theta )\int_{\theta
}^{b}I_{r-}^{2\alpha -1,\theta }\left( \frac{x_{r}-x_{\theta }}{(r-\theta
)^{\alpha +1}}\right) \ D_{b-}^{1-\alpha }y_{b-}(r)drd\theta \\
&=&\frac{\alpha (-1)^{2\alpha -1}}{{\Gamma }(1-\alpha )}%
\int_{a}^{b}D_{a+}^{2\alpha -1}f^{\prime }(x)(\theta )\int_{\theta
}^{b}I_{a+}^{\alpha ,r}I_{r-}^{2\alpha -1,\theta }\left( \frac{%
x_{r}-x_{\theta }}{(r-\theta )^{\alpha +1}}\right) dy_{r}d\theta \,,
\end{eqnarray*}%
where $I^{\alpha ,\theta }$ denotes the fractional integral applied to a
function of $\theta $ and a similar notation is used for fractional
derivatives. Now 
\begin{eqnarray*}
&&I_{a+}^{\alpha ,r}I_{r-}^{2\alpha -1,\theta }\left( \frac{x_{r}-x_{\theta }%
}{(r-\theta )^{\alpha +1}}\right) \\
&=&I_{a+}^{\alpha ,r}\left[ \frac{(-1)^{1-2\alpha }}{{\Gamma }(2\alpha -1)}%
\int_{\theta }^{r}(\theta ^{\prime }-\theta )^{2\alpha -2}\frac{%
x_{r}-x_{\theta ^{\prime }}}{(r-\theta ^{\prime })^{\alpha +1}}d\theta
^{\prime }\right] \\
&=&\frac{(-1)^{1-2\alpha }}{{\Gamma }(2\alpha -1){\Gamma }(\alpha )}%
\int_{\theta <\theta ^{\prime }<\xi <r^{\prime }<r}(r-r^{\prime })^{\alpha
-1}(\theta ^{\prime }-\theta )^{2\alpha -2}(r^{\prime }-\theta ^{\prime
})^{-\alpha -1}dx_{\xi }d\theta ^{\prime }dr^{\prime }\,.
\end{eqnarray*}%
Thus 
\begin{eqnarray}
\lefteqn{A_{2} =\frac{1}{{\Gamma }(2\alpha -1){\Gamma }(\alpha )}%
\int_{a}^{b}D_{a+}^{2\alpha -1}f^{\prime }(x)(\theta ) }    \label{2.12} \\
&&\times \int_{\theta }^{b}\left( \int_{\theta <\theta ^{\prime }<\xi
<r^{\prime }<\eta }(\eta -r^{\prime })^{\alpha -1}(\theta ^{\prime }-\theta
)^{2\alpha -2}(r^{\prime }-\theta ^{\prime })^{-\alpha -1}d\theta ^{\prime
}dr^{\prime }\right) dx_{\xi }dy_{\eta }d\theta \,.  \notag
\end{eqnarray}%
Making the change of variable $\frac{r^{\prime }-\xi }{\eta -\xi }=w$ and
using formula 3.196 in Gradshteyn and Ryzhik \cite{GR} we obtain 
\begin{eqnarray*}
&&\int_{\xi }^{\eta }\left( \eta -r^{\prime }\right) ^{\alpha -1}\left(
r^{\prime }-\theta ^{\prime }\right) ^{-\alpha -1}dr^{\prime } \\
&=&\left( \eta -\xi \right) ^{-1}\int_{0}^{1}(1-w)^{\alpha -1}(w+\frac{\xi
-\theta ^{\prime }}{\eta -\xi })^{-\alpha -1}dw \\
&=&\frac{1}{\alpha }\left( \eta -\xi \right) ^{-1}\left( \frac{\xi -\theta
^{\prime }}{\eta -\xi }\right) ^{-\alpha -1}F(1,\alpha +1,\alpha +1,-\frac{%
\eta -\xi }{\xi -\theta ^{\prime }}) \\
&=&\frac{1}{\alpha }\left( \eta -\xi \right) ^{-1}\left( \frac{\xi -\theta
^{\prime }}{\eta -\xi }\right) ^{-\alpha -1}(1+\frac{\eta -\xi }{\xi -\theta
^{\prime }})^{-1} \\
&=&\frac{1}{\alpha }\left( \eta -\xi \right) ^{\alpha }\left( \xi -\theta
^{\prime }\right) ^{-\alpha }\left( \eta -\theta ^{\prime }\right) ^{-1},
\end{eqnarray*}%
and substituting this expression into (\ref{2.12}) yields%
\begin{eqnarray*}
A_{2} &=&\frac{1}{{\Gamma }(2\alpha -1){\Gamma }(\alpha )\alpha }%
\int_{a}^{b}D_{a+}^{2\alpha -1}f^{\prime }(x)(\theta ) \\
&&\times \int_{\theta <\xi <\eta <b}\left( \eta -\xi \right) ^{\alpha
}\left( \int_{\theta }^{\xi }\left( \xi -\theta ^{\prime }\right) ^{-\alpha
}\left( \eta -\theta ^{\prime }\right) ^{-1}(\theta ^{\prime }-\theta
)^{2\alpha -2}d\theta ^{\prime }\right) dx_{\xi }dy_{\eta }d\theta \,.
\end{eqnarray*}%
Using (\ref{2.6}) we get 
\begin{eqnarray*}
&&\left( \eta -\xi \right) ^{\alpha }\int_{\theta }^{\xi }\left( \xi -\theta
^{\prime }\right) ^{-\alpha }\left( \eta -\theta ^{\prime }\right)
^{-1}(\theta ^{\prime }-\theta )^{2\alpha -2}d\theta ^{\prime } \\
&=&(\eta -\xi )^{\alpha -1}(\xi -\theta )^{\alpha
-1}\int_{0}^{1}(1-q)^{-\alpha }q^{2\alpha -2}\left( 1+(1-q)\frac{\xi -\theta 
}{\eta -\xi }\right) ^{-1}dq \\
&=&G(\theta ,\xi ,\eta ){\Gamma }(2\alpha -1){\Gamma }(\alpha )\alpha.
\end{eqnarray*}%
Hence,%
\begin{equation*}
A_{2}=\int_{a}^{b}D_{a+}^{2\alpha -1}f^{\prime }(x)(\theta )\int_{\theta
<\xi <\eta <b}G(\theta ,\xi ,\eta )dx_{\xi }dy_{\eta }d\theta.
\end{equation*}

Applying the two-dimensional fractional integration by parts formula (\ref%
{1.9}) to $\varphi (\xi ,\eta )=G(\theta ,\xi ,\eta )$ and $\psi (\xi ,\eta
)=(x\otimes y)_{\xi ,\eta }$ and \ using (\ref{2.7}) we obtain%
\begin{eqnarray*}
\int_{\theta< \xi < \eta < b}G(\theta ,\xi ,\eta )dx_{\xi
}dy_{\eta } &=&-\int_{\theta }^{b}\int_{\theta }^{\eta }D_{\theta +}^{\alpha
-{\varepsilon ,\xi }}D_{b-}^{\alpha -{\varepsilon ,\eta }} 
G_{b-}(\theta ,\xi ,\eta )  \\
&&\times \Gamma ^{\alpha -{\varepsilon }}(x\otimes y)_{\xi ,\eta }d\xi d\eta
\\
&=&\int_{\theta }^{b}\int_{\theta }^{\eta }K_{\theta ,b}(\xi ,\eta )\Gamma
^{\alpha -\varepsilon }(x\otimes y)_{\xi ,\eta }d\xi d\eta \,.
\end{eqnarray*}%
Thus 
\begin{equation*}
A_{2}=\int_{a}^{b}D_{a+}^{2\alpha -1}f^{\prime }\left( x\right)
(r)\int_{r}^{b}\int_{r}^{\eta }K_{r,b}(\xi ,\eta )\Gamma ^{\alpha
-\varepsilon }(x\otimes y)_{\xi ,\eta }d\xi d\eta dr.
\end{equation*}%
This proves the theorem.
\end{proof}

For any $(s,t)\in \Delta $, and given a $\beta $-H\"{o}lder continuous
multiplicative functional $(x,y,x\otimes y)$, we define

\begin{eqnarray}
\left\| x\right\| _{s,t,\beta } &=&\sup_{s\leq \theta <r\leq t}\frac{%
|x_{r}-x_{\theta }|}{|r-\theta |^{\beta }},  \label{2.13} \\
\left\| x\otimes y\right\| _{s,t,\beta } &=&\sup_{s\leq \theta <r\leq t}%
\frac{|\left( x\otimes y\right) _{\theta ,r}|}{|r-\theta |^{2\beta }}.
\label{2.14}
\end{eqnarray}%
We also set $\left\| x\right\| _{\beta }=\left\| x\right\| _{0,T,\beta }$,
and $\left\| x\otimes y\right\| _{\beta }=\left\| x\otimes y\right\|
_{0,T,\beta }$. Also, $\left\| \cdot \right\| _{s,t,\infty }$ will denote
the supremum norm in the interval $[s,t]$. In the sequel, $k$ will denote a
constant that may depend on the parameters $\beta $, $\alpha $, $\lambda $, $%
\varepsilon$ and $T$.

The following estimate is useful.

\begin{proposition}
\label{p.3.2} Under the hypotheses of Definition \ref{def1} we have, if $%
b-a\leq 1$ 
\begin{eqnarray}
\left\| \int f(x_{r})dy_{r}\right\| _{a,b,\beta } &\leq &k\left\| f\right\|
_{\infty }\left\| y\right\| _{a,b,\beta }+k\left( \left\| x\otimes y\right\|
_{a,b,\beta }+\left\| x\right\| _{a,b,\beta }\left\| y\right\| _{a,b,\beta
}\right)  \notag \\
&&\!\!\!\!\times \left( \left\| f^{\prime }\right\| _{\infty }+\left\| f^{\prime
}\right\| _{\lambda }\left\| x\right\| _{a,b,\beta }^{\lambda
}(b-a)^{\lambda \beta }\right) (b-a)^{\beta -2\varepsilon }.  \label{2.15}
\end{eqnarray}%
Moreover, if the second derivative $f^{\prime \prime }$ is $\lambda $-H\"{o}%
lder continuous and bounded, and $(\tilde{x},y,\tilde{x}\otimes y)$ is also
a $\beta $-\textit{H\"{o}lder continuous multiplicative functional }on%
\textit{\ }$\mathbb{R}^{m}\otimes \mathbb{R}^{d}$, then%
\begin{eqnarray}
&&\left\| \int f(x_{r})dy_{r}-\int f(\tilde{x}_{r})dy_{r}\right\|
_{a,b,\beta }  \notag \\
&\leq &kH_{1}\Vert x-\tilde{x}\Vert _{a,b,\infty }+kH_{2}\Vert x-\tilde{x}%
\Vert _{a,b,\beta }+kH_{3}\Vert (x-\tilde{x})\otimes y\Vert _{a,b,{\beta }%
},\,  \label{2.16}
\end{eqnarray}%
where%
\begin{eqnarray*}
H_{1} &=&\left\| y\right\| _{a,b,\beta }\left( \Vert f^{\prime }\Vert
_{\infty }+\left\| f^{\prime \prime }\right\| _{\lambda }\Vert \widetilde{x}%
\Vert _{a,b,{\beta }}\left( \Vert x\Vert _{a,b,{\beta }}^{\lambda }+\Vert 
\tilde{x}\Vert _{a,b,{\beta }}^{\lambda }\right) \right) (b-a)^{\beta
(1+\lambda )} \\
&&+\left( \Vert f^{\prime \prime }\Vert _{\infty }+\Vert f^{\prime \prime
}\Vert _{\lambda }\left( \Vert x\Vert _{a,b,\beta }^{\lambda }+\Vert 
\widetilde{x}\Vert _{a,b,\beta }^{\lambda }\right) (b-a)^{\beta \lambda
}\right) \\
&&\times \left( \Vert x\otimes y\Vert _{a,b,{\beta }}+\Vert x\Vert _{a,b,{%
\beta }}\left\| y\right\| _{a,b,\beta }\ \right) (b-a)^{\beta -2\varepsilon
}, \\
H_{2} &=&\ \left\| f^{\prime \prime }\right\| _{\infty }\left\| y\right\|
_{a,b,\beta }\left( \Vert x\Vert _{a,b,{\beta }}+\Vert \tilde{x}\Vert _{a,b,{%
\beta }}\right) (b-a)^{\beta (1+\lambda )} \\
&&+\left\| f^{\prime \prime }\right\| _{\infty }\left( \Vert x\otimes y\Vert
_{a,b,{\beta }}+\Vert x\Vert _{a,b,{\beta }}\left\| y\right\| _{a,b,\beta }\
\right) (b-a)^{2\beta -2\varepsilon } \\
&&+\ \left( \Vert f\Vert _{\infty }+\Vert f\Vert _{\lambda }\Vert \tilde{x}%
\Vert _{a,b,{\beta }}^{\lambda }(b-a)^{\lambda \beta }\right) \left\|
y\right\| _{a,b,\beta }(b-a)^{\beta -2\varepsilon } \\
H_{3} &=&\left( \Vert f\Vert _{\infty }+\Vert f\Vert _{\lambda }\Vert \tilde{%
x}\Vert _{a,b,{\beta }}^{\lambda }(b-a)^{\lambda \beta }\right) (b-a)^{\beta
-2\varepsilon }.
\end{eqnarray*}
\end{proposition}

\bigskip

\textbf{Remark: }In (\ref{2.15}) we can replace $\left\| f\right\| _{\infty
} $ and $\left\| f^{\prime }\right\| _{\infty }$ by $\left\| f(x)\right\|
_{a,b,\infty }$ and $\left\| f^{\prime }(x)\right\| _{a,b,\infty }$,
respectively.

\medskip

\begin{proof}
First we have, for any $r\in \lbrack a,b]$ 
\begin{equation}
\left| D_{a+}^{\alpha }f(x)(r)\right| \leq k\left( |f(x_{r})|\left(
r-a\right) ^{-\alpha }+\left\| f\right\| _{\lambda }\left\| x\right\|
_{a,r,\beta }^{\lambda }\left( r-a\right) ^{\lambda \beta -\alpha }\right)
\,,  \label{2.17}
\end{equation}%
\begin{equation}
\left| \hat{D}_{a+}^{\alpha }f(x)(r)\right| \leq k\left( |f(x_{r})|\left(
r-a\right) ^{-\alpha }+\left\| f^{\prime }\right\| _{\lambda }\left\|
x\right\| _{a,r,\beta }^{1+\lambda }\left( r-a\right) ^{(\lambda +1)\beta
-\alpha }\right) ,  \label{2.18}
\end{equation}%
and 
\begin{equation}
\left| D_{b-}^{1-\alpha }y_{b-}(r)\right| \leq k\left\| y\right\|
_{r,b,\beta }(b-r)^{\alpha +\beta -1}\,.  \label{2.19}
\end{equation}%
The expression (\ref{1.11}) of $\Gamma $ yields 
\begin{equation}
\left| \Gamma ^{\alpha -\varepsilon }\left( x\otimes y\right) _{a,b}\right|
\leq k\left( \left\| x\otimes y\right\| _{a,b,\beta }+\left\| x\right\|
_{a,b,\beta }\left\| y\right\| _{a,b,\beta }\right) (b-a)^{2\beta +2\alpha
-2-2\varepsilon }.  \label{2.20}
\end{equation}%
Consequently, from (\ref{2.8}) Lemma \ref{l.7.2} we obtain the estimate 
\begin{equation}
|{\Lambda }_{a}^{b}(x\otimes y)|\leq k\left( \left\| x\otimes y\right\|
_{a,b,\beta }+\left\| x\right\| _{a,b,\beta }\left\| y\right\| _{a,b,\beta
}\right) (b-a)^{2\beta +2\alpha -2-2\varepsilon }.  \label{2.21}
\end{equation}%
Thus 
\begin{eqnarray*}
\lefteqn{\left| \int_{a}^{b}f(x_{r})dy_{r}\right| \leq k\left\| y\right\|
_{a,b,\beta }\ \left( \int_{a}^{b}|f(x_{r})|\left( r-a\right) ^{-\alpha
}(b-r)^{\alpha +\beta -1}dr\right.} \\
&&\qquad \left. +\left\| f^{\prime }\right\| _{\lambda }\left\| x\right\|
_{a,b,\beta }^{1+\lambda }\int_{a}^{b}\left( r-a\right) ^{(\lambda +1)\beta
-\alpha }(b-r)^{\alpha +\beta -1}dr\right) \\
&&\qquad +k\int_{a}^{b}\left( \frac{\left\| f^{\prime }\right\| _{\infty }}{%
(r-a)^{2\alpha -1}}+\left\| f^{\prime }\right\| _{\lambda }\left\| x\right\|
_{a,b,\beta }^{\lambda }\left( r-a\right) ^{\lambda \beta -2\alpha +1}\right)
\\
&&\qquad \quad \times \ \left( \left\| x\otimes y\right\| _{a,b,\beta
}+\left\| x\right\| _{a,b,\beta }\left\| y\right\| _{a,b,\beta }\right)
(b-r)^{2\beta +2\alpha -2-2\varepsilon }dr.
\end{eqnarray*}%
Therefore, we obtain%
\begin{eqnarray*}
\lefteqn{\left| \int_{a}^{b}f(x_{r})dy_{r}\right| \leq k\left\| f\right\| _{\infty
}\left\| y\right\| _{a,b,\beta }(b-a)^{\beta } }\\
&&+k\left( \left\| x\otimes y\right\| _{a,b,\beta }+\left\| x\right\|
_{a,b,\beta }\left\| y\right\| _{a,b,\beta }\right) \left\| f^{\prime
}\right\| _{\infty }(b-a)^{2\beta -2\varepsilon } \\
&&+k\left\| f^{\prime }\right\| _{\lambda }\left\| x\right\| _{a,b,\beta
}^{\lambda }(\left\| x\right\| _{a,b,\beta }\left\| y\right\| _{a,b,\beta
}+\left\| x\otimes y\right\| _{a,b,\beta })(b-a)^{\left( \lambda +2\right)
\beta -2\varepsilon }
\end{eqnarray*}%
and this implies (\ref{2.15}) easily.

Note that for any $a\leq \theta \leq r\leq b$ we can write 
\begin{eqnarray*}
&&f(x_{r})-f(x_{\theta })-f^{\prime }(x_{\theta })(x_{r}-x_{\theta })-\left[
f(\tilde{x}_{r})-f(\tilde{x}_{\theta })-f^{\prime }(\tilde{x}_{\theta })(%
\tilde{x}_{r}-\tilde{x}_{\theta })\right] \\
&=&\int_{0}^{1}\left[ f^{\prime }(x_{\theta }+z(x_{r}-x_{\theta
}))-f^{\prime }(x_{\theta })\right] (x_{r}-x_{\theta }-\widetilde{x}_{r}+%
\widetilde{x}_{\theta })dz \\
&&+\int_{0}^{1}\left[ f^{\prime }(x_{\theta }+z(x_{r}-x_{\theta
}))-f^{\prime }(\widetilde{x}_{\theta }+z(\widetilde{x}_{r}-\widetilde{x}%
_{\theta }))-f^{\prime }(x_{\theta })+f^{\prime }(\widetilde{x}_{\theta })%
\right] (\widetilde{x}_{r}-\widetilde{x}_{\theta })dz \\
&=&a_{1}-a_{2}.
\end{eqnarray*}%
We have%
\begin{equation*}
\left| a_{1}\right| \leq \left\| f^{\prime \prime }\right\| _{\infty }\Vert
x\Vert _{a,b,{\beta }}\Vert x-\tilde{x}\Vert _{a,b,{\beta }}(r-\theta )^{2{%
\beta }}.
\end{equation*}%
For the term $a_{2}$ we make the decomposition%
\begin{eqnarray*}
\lefteqn{a_{2} =(\widetilde{x}_{r}-\widetilde{x}_{\theta })\int_{0}^{1}\int_{0}^{1}%
\left[ \left( \widetilde{x}_{\theta }-x_{\theta }\right) +z\left( \widetilde{%
x}_{r}-\widetilde{x}_{\theta }-x_{r}+x_{\theta }\right) \right] }\\
&&\qquad \times f^{\prime \prime }(x_{\theta }+z(x_{r}-x_{\theta })+t(\widetilde{x}%
_{\theta }-x_{\theta })+tz(\widetilde{x}_{r}-\widetilde{x}_{\theta
}-x_{r}+x_{\theta }))dtdz \\
&&\qquad -(\widetilde{x}_{r}-\widetilde{x}_{\theta })(\widetilde{x}_{\theta
}-x_{\theta })\int_{0}^{1}f^{\prime \prime }(x_{\theta }+t(\widetilde{x}%
_{\theta }-x_{\theta }))dt\\
&&= \left( \widetilde{%
x}_{r}-\widetilde{x}_{\theta }-x_{r}+x_{\theta }\right)    \\
&&\quad \times
\int_{0}^{1}\int_{0}^{1}z%
 f^{\prime \prime }(x_{\theta }+z(x_{r}-x_{\theta })+t(\widetilde{x}%
_{\theta }-x_{\theta })+tz(\widetilde{x}_{r}-\widetilde{x}_{\theta
}-x_{r}+x_{\theta }))dtdz \\
&&-(\widetilde{x}_{r}-\widetilde{x}_{\theta })(\widetilde{x}_{\theta
}-x_{\theta }) 
 \int_{0}^{1}\int_{0}^{1}%
\left[   f^{\prime \prime }(x_{\theta }+z(x_{r}-x_{\theta })+t(\widetilde{x}%
_{\theta }-x_{\theta })+tz(\widetilde{x}_{r}-\widetilde{x}_{\theta
}-x_{r}+x_{\theta })) \right. \\
&&\quad \left. - f^{\prime \prime }(x_{\theta }+t(\widetilde{x}%
_{\theta }-x_{\theta }))\right] dt dz.
\end{eqnarray*}%
Thus,%
\begin{eqnarray*}
\left| a_{2}\right| &\leq &\left\| f^{\prime \prime }\right\| _{\infty
}\Vert \widetilde{x}\Vert _{a,b,{\beta }}\Vert x-\tilde{x}\Vert _{a,b,{\beta 
}}(r-\theta )^{2{\beta }} \\
&&+\Vert f^{\prime \prime }\Vert _{\lambda }\Vert \tilde{x}\Vert _{a,b,{%
\beta }}\Vert x-\tilde{x}\Vert _{a,b,{\infty }}\left( \Vert x\Vert _{a,b,{%
\beta }}^{\lambda }+\Vert x-\tilde{x}\Vert _{a,b,{\beta }}^{\lambda }\right)
(r-\theta )^{{\beta (1+\lambda )}}.
\end{eqnarray*}%
As a consequence,%
\begin{eqnarray}
&&\left| f(x_{r})-f(x_{\theta })-f^{\prime }(x_{\theta })(x_{r}-x_{\theta
})- \left[ f(\tilde{x}_{r})-f(\tilde{x}_{\theta })-f^{\prime }(\tilde{x}%
_{\theta })(\tilde{x}_{r}-\tilde{x}_{\theta })\right] \right|  \notag \\
&\leq &kI_{1}(r-\theta )^{{\beta (1+\lambda )}},  \label{2.22}
\end{eqnarray}%
where 
\begin{eqnarray*}
I_{1} &=&\ \Vert f^{\prime \prime }\Vert _{\infty }\left\{ \Vert x\Vert
_{a,b,{\beta }}+\Vert \tilde{x}\Vert _{a,b,{\beta }}\right\} \Vert x-\tilde{x%
}\Vert _{a,b,{\beta }} \\
&&+\Vert f^{\prime \prime }\Vert _{\lambda }\Vert \widetilde{x}\Vert _{a,b,{%
\beta }}\left( \Vert x\Vert _{a,b,{\beta }}^{\lambda }+\Vert \tilde{x}\Vert _{a,b,{%
\beta }}^{\lambda } \right)  \Vert x-\tilde{x}\Vert _{a,b,\infty }\,.
\end{eqnarray*}%
On the other hand, we have 
\begin{eqnarray*}
&&D_{s+}^{2\alpha -1}f^{\prime }(x)(r)-D_{s+}^{2\alpha -1}f^{\prime }(\tilde{%
x})(r) \\
&=&\frac{1}{\Gamma (2-2\alpha )}\left\{ \frac{f^{\prime }(x_{r})-f^{\prime }(%
\tilde{x}_{r})}{(r-s)^{2\alpha -1}}\right. \\
&&\left. +(2\alpha -1)\int_{s}^{r}\frac{\left[ f^{\prime }(x_{r})-f^{\prime
}(\widetilde{x}_{r})-f^{\prime }(x_{\theta })+f^{\prime }(\widetilde{x}%
_{\theta })\right] }{(r-\theta )^{2\alpha }}d\theta \right\} .
\end{eqnarray*}%
Using the decomposition%
\begin{eqnarray*}
&&f^{\prime }(x_{r})-f^{\prime }(\widetilde{x}_{r})-f^{\prime }(x_{\theta
})+f^{\prime }(\widetilde{x}_{\theta }) \\
&=&\int_{0}^{1}f^{\prime \prime }(x_{r}+t(\widetilde{x}_{r}-x_{r}))(%
\widetilde{x}_{r}-x_{r})dt-\int_{0}^{1}f^{\prime \prime }(x_{\theta }+t(%
\widetilde{x}_{\theta }-x_{\theta }))(\widetilde{x}_{\theta }-x_{\theta })dt,
\end{eqnarray*}%
we obtain 
\begin{eqnarray}
&&\left| D_{s+}^{2\alpha -1}f^{\prime }(x)(r)-D_{s+}^{2\alpha -1}f^{\prime }(%
\tilde{x})(r)\right|  \notag \\
&\leq &k(r-s)^{1-2\alpha }\Vert f^{\prime \prime }\Vert _{\infty }\Vert x-%
\tilde{x}\Vert _{s,r,\infty }+k(r-s)^{\beta -2\alpha +1}\Vert f^{\prime
\prime }\Vert _{\infty }\Vert x-\tilde{x}\Vert _{s,r,\beta }  \notag \\
&&\quad +k(r-s)^{\beta \lambda -2\alpha +1}\Vert f^{\prime \prime }\Vert
_{\lambda }\left( \Vert x\Vert _{s,r,\beta }^{\lambda }+\Vert \widetilde{x}\Vert
_{s,r,\beta }^{\lambda }\right)  \Vert x-\tilde{x}\Vert _{s,r,\infty }  \notag
\\
&=&kI_{2}(r-s)^{1-2\alpha },  \label{2.23}
\end{eqnarray}%
where 
\begin{eqnarray*}
I_{2} &=&\left( \Vert f^{\prime \prime }\Vert _{\infty }+\Vert f^{\prime
\prime }\Vert _{\lambda }\left( \Vert x\Vert _{a,b,\beta }^{\lambda }+\Vert \widetilde{x%
}\Vert _{a,b,\beta }^{\lambda }\right)  (b-a)^{\beta \lambda }\right) \Vert x-%
\tilde{x}\Vert _{a,b,\infty } \\
&&+\Vert f^{\prime \prime }\Vert _{\infty }\Vert x-\tilde{x}\Vert
_{a,b,\beta }(b-a)^{\beta }.
\end{eqnarray*}%
Now using (\ref{2.19}), (\ref{2.22}), (\ref{2.23}) we obtain 
\begin{eqnarray}
\lefteqn{\left| \int_{a}^{b}\left[ f(x_{r})-f(\tilde{x}_{r})\right] dy_{r}\right|
\leq k\left\| y\right\| _{a,b,\beta}   } \notag\\
&&\times \left( \int_{a}^{b}|f(x_{r})-f(%
\tilde{x}_{r})|\left( r-a\right) ^{-\alpha }(b-r)^{\alpha +\beta -1}dr\right.
\notag \\
&& \left. +I_{1}\int_{a}^{b}\left( r-a\right) ^{\beta (1+\lambda
)-\alpha }(b-r)^{\alpha +\beta -1}dr\right)  \notag \\
&& +kI_{2}\int_{a}^{b}(b-r)^{1-2\alpha }\left| \Lambda
_{r}^{b}(x\otimes y)\right| dr  \notag \\
&& +\int_{a}^{b}\left| D_{a+}^{2\alpha -1}f^{\prime }\left( \tilde{x}%
\right) (r)\right| \left| \Lambda _{r}^{b}(\left[ x-\tilde{x}\right] \otimes
y)\right| dr.  \notag
\end{eqnarray}%
Finally, using (\ref{2.21}) we get 
\begin{eqnarray}
\lefteqn{\left| \int_{a}^{b}\left[ f(x_{r})-f(\tilde{x}_{r})\right] dy_{r}\right|
\leq k\left\| y\right\| _{a,b,\beta } } \notag \\
&&\left( \Vert f^{\prime }\Vert
_{\infty }(b-a)^{\beta }\Vert x-\tilde{x}\Vert _{a,b,\infty }+I_{1}(b-a)^{{%
\beta (2+\lambda )}}\right)  \notag \\
&&+kI_{2}\left( \Vert x\otimes y\Vert _{a,b,{\beta }}+\Vert x\Vert _{a,b,{%
\beta }}\Vert y\Vert _{a,b,{\beta }}\right) (b-a)^{2{\beta -2\varepsilon }} 
\notag \\
&&+\left[ \Vert f\Vert _{\infty }+\Vert f\Vert _{\lambda }\Vert \tilde{x}%
\Vert _{a,b,{\beta }}^{\lambda }(b-a)^{{\lambda }{\beta }}\right]
\label{2.24} \\
&&\times \left\{ \Vert (x-\tilde{x})\otimes y\Vert _{a,b,{\beta }}+\Vert x-%
\tilde{x}\Vert _{a,b,{\beta }}\Vert y\Vert _{a,b,{\beta }}\right\} (b-a)^{2{%
\beta -2\varepsilon }}\,.  \notag
\end{eqnarray}%
This implies (\ref{2.16}).
\end{proof}

The following corollary is the direct consequence of the proposition.

\begin{corollary}
\label{c.3.3} Assume $b-a\leq 1$. Under the hypotheses of Definition \ref%
{def1}, if $(x,\tilde{y},x\otimes \tilde{y})$ is also a $\beta $-\textit{H%
\"{o}lder continuous multiplicative functional }on\textit{\ }$\mathbb{R}%
^{m}\otimes \mathbb{R}^{d}$, we have%
\begin{eqnarray}
&&\left\| \int f(x_{r})dy_{r}-\int f(x_{r})d\tilde{y}_{r}\right\|
_{a,b,\beta }  \notag \\
&\leq &k\left\| f\right\| _{\infty }\left\| y-\tilde{y}\right\| _{a,b,\beta
}+k\left( \left\| x\otimes (y-\tilde{y})\right\| _{a,b,\beta }+\left\|
x\right\| _{a,b,\beta }\left\| y-\tilde{y}\right\| _{a,b,\beta }\right) 
\notag \\
&&\times \left( \left\| f^{\prime }\right\| _{\infty }+\left\| f^{\prime
}\right\| _{\lambda }\left\| x\right\| _{a,b,\beta }^{\lambda
}(b-a)^{\lambda \beta }\right) (b-a)^{\beta -2\varepsilon }.  \label{2.25}
\end{eqnarray}%
On the other hand, if the derivative $f^{\prime \prime }$ is $\lambda $-H%
\"{o}lder continuous and bounded, $(\tilde{x},y,\tilde{x}\otimes y)$ is
another $\beta $-\textit{H\"{o}lder continuous multiplicative functional }on%
\textit{\ }$\mathbb{R}^{m}\otimes \mathbb{R}^{d}$, and $\tilde{f}$ is
another function satisfying the hypotheses of Definition \ref{def1}, then 
\begin{eqnarray}
&&\left\| \int f(x_{r})dy_{r}-\int \tilde{f}(\tilde{x}_{r})dy_{r}\right\|
_{a,b,\beta }  \notag \\
&\leq &kH_{1}^{f}\Vert x-\tilde{x}\Vert _{a,b,\infty }+kH_{2}^{f}\Vert x-%
\tilde{x}\Vert _{a,b,\beta }+kH_{3}^{f}\Vert (x-\tilde{x})\otimes y\Vert
_{a,b,{\beta }}  \notag \\
&&+k\left\| f-\widetilde{f}\right\| _{\infty }\left\| y\right\| _{a,b,\beta
}+k\left( \left\| x\otimes y\right\| _{a,b,\beta }+\left\| x\right\|
_{a,b,\beta }\left\| y\right\| _{a,b,\beta }\right)  \notag \\
&&\times \left( \left\| f^{\prime }-\tilde{f}^{\prime }\right\| _{\infty
}+\left\| f^{\prime }-\tilde{f}^{\prime }\right\| _{\lambda }\left\| \tilde{x%
}\right\| _{a,b,\beta }^{\lambda }(b-a)^{\lambda \beta }\right) (b-a)^{\beta
-2\varepsilon }.  \label{e.3.29}
\end{eqnarray}
\end{corollary}

The estimate (\ref{2.25}) implies that for a fixed $x$, the mapping $%
(y,x\otimes y)\rightarrow \int f(x_{r})dy_{r}$ is continuous with respect to
the $\beta $-norm. As a consequence, if $y^{n}$ is a sequence of
continuously differentiable functions (or Lipschitz functions) such that%
\begin{eqnarray*}
\left\| y-y^{n}\right\| _{\beta } &\rightarrow &0, \\
\left\| x\otimes y-x\otimes y^{n}\right\| _{\beta } &\rightarrow &0
\end{eqnarray*}%
as $n$ tends to infinity, then%
\begin{equation}
\left\| \int f(x_{r})dy_{r}-\int f(x_{r})dy_{r}^{n}\right\| _{\beta
}\rightarrow 0.  \label{2.26}
\end{equation}%
Hence, the integral $\int f(x_{r})dy_{r}$ introduced in Definition \ref{def1}
does not depend on the parameters $\alpha $ and $\varepsilon $, and it
coincides with the classical integral $\int f(x_{r})y_{r}^{\prime }dr$ when $%
y$ is continuously differentiable.

Set $t_{i}^{n}=\frac{iT}{n}$ for $i=0,1,\ldots ,n$. If $y$ is $\beta $-H\"{o}%
lder continuous, the sequence of functions 
\begin{equation*}
y_{t}^{n}=y_{0}\mathbf{1}_{\{0\}}(t)+\sum_{i=1}^{n}\mathbf{1}%
_{(t_{i-1}^{n},t_{i}^{n}]}(t)\left[ y_{t_{i-1}^{n}}+\frac nT\left(
t-t_{i-1}^{n}\right) (y_{t_{i}^{n}}-y_{t_{i-1}^{n}})\right]
\end{equation*}%
converge to $y$ in the $\beta' $-norm for any $\beta'<\beta$. Assume that the multiplicative
functional $\int_{s}^{t}(x_{r}-x_{s})dy_{r}^{n}$ converges in the $\beta' $%
-norm as $n$ tends to infinity to $\left( x\otimes y\right) _{s,t}$. Then (%
\ref{2.26}) holds with $\beta=\beta'$. In particular, this means that 
\begin{equation}
\int_{0}^{T}f(x_{r})dy_{r}=\lim_{n\rightarrow \infty }\sum_{i=1}^{n}\frac{n}{%
T}\left( \int_{t_{i-1}^{n}}^{t_{i}^{n}}f(x_{s})ds\right)
(y_{t_{i}^{n}}-y_{t_{i-1}^{n}}).  \label{eqa}
\end{equation}

For any $p\geq 1$, the $p$ variation of a function $x:[0,T]\rightarrow 
\mathbb{R}$ is defined as%
\begin{equation*}
\mathrm{Var}_{p}(x)=\sup_{\pi }\left( \sum_{i=1}^{n}\left\vert
x(t_{i}^{n})-x(t_{i-1}^{n})\right\vert ^{p}\right) ^{1/p},
\end{equation*}%
where $\pi =\{0=t_{0}<\cdots <t_{n}=T\}$ runs over all partitions of $[0,T]$%
. Notice that 
\begin{equation*}
\mathrm{Var}_{1/\beta }(x)\leq \left\Vert x\right\Vert _{\beta }.
\end{equation*}%
Then, for any $\beta $-H\"{o}lder continuous multiplicative functional $%
(x,y,x\otimes y)$ on $\mathbb{R}^{m}\otimes \mathbb{R}^{d}$ and any function 
$f$ satisfying the hypotheses of Definition \ref{def1}, the integral $%
\int_{0}^{T}f(x_{r})dy_{r}$ coincides with the integral defined using the $%
\frac{1}{\beta }$-variation norm (see \cite{LQ}). This implies that $%
\int_{0}^{T}f(x_{s})dy_{s}$ is given by the limit of the Riemann sums of the
form%
\begin{equation*}
\int_{0}^{T}f(x_{s})dy_{s}=\lim_{\left\vert \pi \right\vert \rightarrow
0}\sum_{i=1}^{n}f(x_{t_{i-1}})(y_{t_{i}}-y_{t_{i-1}})+f^{\prime
}(x_{t_{i-1}})(x\otimes y)_{t_{i-1},t_{i}},
\end{equation*}%
where $\pi =\{0=t_{0}<\cdots <t_{n}=T\}$ runs over all partitions of $[0,T]$.

\bigskip

In order to handle differential equations we need to introduce the tensor
product of two multiplicative functionals:

\begin{definition}
\label{d.3.3} Suppose that $(x,y,x\otimes y)$ and $(y,z,y\otimes z)$ are $%
\beta $-H\"{o}lder continuous real valued multiplicative functionals. Then,
for all $\ a\leq b\leq c$, we define 
\begin{eqnarray*}
\lefteqn{\left( x\otimes \left( y\otimes z\right) _{\cdot ,c}\right) _{a,b}
=\int_{a}^{b}\Lambda _{a,r,b,c}(x,z)D_{b-}^{1-\alpha }y_{b-}(r)dr} \\
&&+\frac{1}{\Gamma (2-2\alpha )}\int_{a<r<\xi <\eta <b}K_{r,b}(\xi ,\eta ) \\
&&\times \frac{\Gamma ^{\alpha -\varepsilon }\left( x\otimes y\right) _{\xi
,\eta }\left( z_{c}-z_{r}\right) -\left( x_{r}-x_{a}\right) \Gamma ^{\alpha
-\varepsilon }\left( y\otimes z\right) _{\xi ,\eta }}{(r-a)^{2\alpha -1}}%
drd\xi d\eta ,
\end{eqnarray*}%
where%
\begin{equation*}
\Lambda _{a,r,b,c}(x,z)=\frac{(-1)^{\alpha }}{\Gamma \left( 1-\alpha \right) 
}\left( \frac{(x_{r}-x_{a})(z_{c}-z_{r})}{(r-a)^{\alpha }}+\alpha
\int_{a}^{r}\frac{(z_{r}-z_{\theta })(x_{\theta }-x_{r})}{\left( r-\theta
\right) ^{\alpha +1}}d\theta \right) .
\end{equation*}
\end{definition}

We have the following result.

\begin{proposition}
\label{p.3.4} If the function $y$ is continuously differentiable and for all 
$a\leq b$ 
\begin{eqnarray*}
\left( y\otimes z\right) _{a,b} &=&\int_{a}^{b}(z_{b}-z_{r})\ y_{r}^{\prime
}dr, \\
\left( x\otimes y\right) _{a,b} &=&\int_{a}^{b}\ (x_{r}-x_{a})y_{r}^{\prime
}dr,
\end{eqnarray*}%
then%
\begin{equation*}
\left( x\otimes \left( y\otimes z\right) _{\cdot ,c}\right)
_{a,b}=\int_{a}^{b}(x_{r}-x_{a})(z_{c}-z_{r})y_{r}^{\prime }dr.
\end{equation*}
\end{proposition}

\begin{proof}
We are going to use formula (\ref{2.9}) with $m=2$, $d=1$, $f(x,z)=xz$
and the functions $x_{t}-x_{a}$ and $z_{c}-z_{t}$. In this way we obtain 
\begin{eqnarray*}
&&\int_{a}^{b}(x_{\theta }-x_{a})(z_{c}-z_{\theta })dy_{\theta } \\
&=&\frac{(-1)^{\alpha }}{\Gamma \left( 1-\alpha \right) }\int_{a}^{b}\left( 
\frac{(x_{r}-x_{a})(z_{c}-z_{r})}{(r-a)^{\alpha }}+\alpha \int_{a}^{r}\frac{%
(z_{r}-z_{\theta })(x_{\theta }-x_{r})}{\left( r-\theta \right) ^{\alpha +1}}%
d\theta \right)     \\
&&\times D_{b-}^{1-\alpha }y_{b-}(r)dr \\
&&+\frac{1}{\Gamma (2-2\alpha )}\int_{a}^{b}\left[ \frac{z_{c}-z_{r}}{%
(r-a)^{2\alpha -1}}+(2\alpha -1)\int_{a}^{r}\frac{z_{r}-z_{\theta }}{%
(r-\theta )^{2\alpha }}d\theta \right] \\
&&\qquad \times \int_{r}^{b}\int_{r}^{\eta }K_{r,b}(\xi ,\eta )\Gamma
^{\alpha -\varepsilon }\left( x\otimes y\right) _{\xi ,\eta }d\xi d\eta dr \\
&&-\frac{1}{\Gamma (2-2\alpha )}\int_{a}^{b}\left[ \frac{x_{r}-x_{a}}{%
(r-a)^{2\alpha -1}}+(2\alpha -1)\int_{a}^{r}\frac{x_{r}-x_{\theta }}{%
(r-\theta )^{2\alpha }}d\theta \right] \\
&&\qquad \times \int_{r}^{b}\int_{r}^{\eta }K_{r,b}(\xi ,\eta )\Gamma
^{\alpha -\varepsilon }\left( y\otimes z\right) _{\xi ,\eta }d\xi d\eta dr,
\end{eqnarray*}%
and this completes the proof.
\end{proof}

It is easy to obtain the following estimate

\begin{proposition}
\label{p.3.8} Suppose that $(x,y,x\otimes y)$ and $(y,z,y\otimes z)$ are $%
\beta $-H\"{o}lder continuous real valued multiplicative functionals. Then,
for any $a\leq b\leq c$ we have 
\begin{eqnarray*}
\lefteqn{\left| \left( x\otimes \left( y\otimes z\right) _{\cdot ,c}\right)
_{a,b}\right|  \leq k\left( \left\| y\right\| _{a,b,\beta }\left\|
x\right\| _{a,b,\beta }\left\| z\right\| _{a,b,\beta }\right.  }\\
&&+\left. \left\| z\right\| _{a,b,\beta }\left\| y\otimes x\right\|
_{a,b,\beta }+\left\| x\right\| _{a,b,\beta }\left\| y\otimes z\right\|
_{a,b,\beta }\right) (b-a)^{3\beta } \\
&&+k\left\| z\right\| _{a,c,\beta }\left( \left\| y\right\| _{a,b,\beta
}\left\| x\right\| _{a,b,\beta }+\ \left\| y\otimes x\right\| _{a,b,\beta }\
\right) (b-a)^{2\beta }(c-a)^{\beta }.
\end{eqnarray*}%
\newline
\end{proposition}

If $b=c$ we write $\left( x\otimes \left( y\otimes z\right) _{\cdot
,b}\right) _{a,b}=\left( x\otimes y\otimes z\right) _{a,b}$. If the
functions $x$, $y$ and $z$ are continuously differentiable, then%
\begin{equation*}
\left( x\otimes y\otimes z\right) _{a,b}=\int_{a<r<\theta <\sigma
<b}x_{t}^{\prime }y_{\theta }^{\prime }z_{\sigma }^{\prime }drd\theta
d\sigma .
\end{equation*}%
Define 
\begin{equation*}
\left\| x\otimes y\otimes z\right\| _{a,b,\beta }=\sup_{a\leq \theta <r\leq
b}\frac{|\left( x\otimes y\otimes z\right) _{\theta ,c}|}{|r-\theta
|^{3\beta }}.
\end{equation*}%
Then, Proposition \ref{p.3.8} implies that 
\begin{eqnarray}
\left\| x\otimes y\otimes z\right\| _{a,b,\beta } &\leq &k\left( \left\|
y\right\| _{a,b,\beta }\left\| x\right\| _{a,b,\beta }\left\| z\right\|
_{a,b,\beta }+\left\| z\right\| _{a,b,\beta }\left\| y\otimes x\right\|
_{a,b,\beta }\right.  \notag \\
&&\qquad +\left. \left\| x\right\| _{a,b,\beta }\left\| y\otimes z\right\|
_{a,b,\beta }\right) .  \label{2.27}
\end{eqnarray}

Proposition \ref{p.3.8} also implies that $(x,\left( y\otimes z\right)
_{\cdot ,c},\left( x\otimes \left( y\otimes z\right) _{\cdot ,c}\right) )$
is a $\beta $-H\"{o}lder \ continuous functional on the interval $[0,c]$. As
a consequence, if $f$ satisfies the assumptions of Definition \ref{def1}, we
can define the integral $\int_{a}^{b}f(x_{r})d_{r}\left( y\otimes z\right)
_{r,c}$, for all $a\leq b\leq c$. The following estimate for this integral
will be needed to solve differential equations.

\begin{proposition}
\label{p.3.1} Suppose that $(x,y,x\otimes y)$ and $(y,z,y\otimes z)$ are $%
\beta $-H\"{o}lder continuous multiplicative functionals on\textit{\ }$%
\mathbb{R}^{m}\otimes \mathbb{R}^{d}$. Let $f:\mathbb{R}^{m}\mathbb{%
\rightarrow R}^{d}$ be a continuously differentiable function such that $%
f^{\prime }$ is $\lambda $-H\"{o}lder continuous, where $\lambda >\frac{1}{%
\beta }-2$. Fix $\alpha >0$ and $\varepsilon >0$ such that $1-\beta <\alpha
<2\beta $, $\alpha <\frac{\lambda \beta +1}{2}$ and $\ $ $\varepsilon
<\alpha +\beta -1$. Then the following estimate holds%
\begin{eqnarray}
&&\sup_{a\leq \xi \leq \eta \leq b}\frac{1}{(\eta -\xi )^{2\beta }}\left|
\int_{\xi }^{\eta }f(x_{r})d_{r}(y\otimes z)_{r,\eta }\right|  \leq \ k\Big[
A_{a,b}   \notag  \\
&& \qquad \qquad+B_{a,b}\left\| x\otimes y\right\| _{a,b,\beta }\,(b-a)^{\beta
-2\varepsilon }\Big] ,  \label{2.28}
\end{eqnarray}%
where 
\begin{eqnarray}
A_{a,b}&= &\left( \Vert y\otimes z\Vert _{a,b,\beta }+\Vert y\Vert
_{a,b,\beta }^{\ }\Vert z\Vert _{a,b,\beta }^{\ }\right)   \label{2.29} \\
&&\times \left[ \left\| f\right\| _{\infty }+\left( \left\| x\right\|
_{a,b,\beta }\left\| f^{\prime }\right\| _{\infty }+\left\| f^{\prime
}\right\| _{\lambda }\left\| x\right\| _{a,b,\beta }^{1+\lambda
}(b-a)^{\lambda \beta }\right) (b-a)^{\beta -2\varepsilon }\right] \,, 
\notag
\end{eqnarray}%
and 
\begin{equation}
B_{a,b}=\ \left\| z\right\| _{a,b,\beta }\left( \left\| f^{\prime }\right\|
_{\infty }+\left\| f^{\prime }\right\| _{\lambda }\left\| x\right\|
_{a,b,\beta }^{\lambda }(b-a)^{\lambda \beta }\right) \,.  \label{2.30}
\end{equation}
\end{proposition}

\begin{proof}
To simplify the proof we will assume $d=m=1$. From (\ref{2.1}) it is easy to
see that 
\begin{equation}
\left\| \left( x\otimes y\right) _{\cdot ,b}\right\| _{a,b,\beta }\leq
\left( \left\| x\otimes y\right\| _{a,b,\beta }+\left\| x\right\|
_{a,b,\beta }\left\| y\right\| _{a,b,\beta }\right) (b-a)^{\beta },
\label{2.31}
\end{equation}%
and from Proposition \ref{p.3.8} we have 
\begin{eqnarray}
&&\Vert x\otimes \left( y\otimes z\right) _{\cdot ,b}\Vert _{a,b,\beta } 
\notag 
\leq k\Big( \left\| x\right\| _{a,b,\beta }\left\| y\right\| _{a,b,\beta
}^{\ }\Vert z\Vert _{a,b,\beta }^{\ }  \\
&&\qquad +\left\| x\right\| _{a,b,\beta }\left\|
y\otimes z\right\| _{a,b,\beta }+\left\| z\right\| _{a,b,\beta }\left\|
x\otimes y\right\| _{a,b,\beta }\Big) (b-a)^{{\beta }}\,.  \label{2.32}
\end{eqnarray}%
From (26), (\ref{2.31}), and (\ref{2.32}) we obtain%
\begin{eqnarray*}
&&\left| \int_{a}^{b}f(x_{r})d_{r}(y\otimes z)_{r,t}\right|  \\
&\leq &k\Vert f\Vert _{\infty }\Vert (y\otimes z)_{\cdot ,b}\Vert
_{a,b,\beta }(b-a)^{\beta } \\
&&+k[\left( \Vert x\otimes (y\otimes z)_{\cdot ,b}\Vert _{a,b,\beta }+\Vert
x\Vert _{a,b,\beta }\Vert (y\otimes z)_{\cdot ,b}\Vert _{a,b,\beta }\right) 
\\
&&\times \left( \Vert f^{\prime }\Vert _{\infty }+\Vert f^{\prime }\Vert _{{%
\lambda }}\Vert x\Vert _{a,b,\beta }^{{\lambda }}(b-a)^{\lambda \beta
}\right) ](b-a)^{2\beta -2\varepsilon } \\
&\leq &k\Vert f\Vert _{\infty }\left( \Vert y\otimes z\Vert _{a,b,\beta
}+\Vert y\Vert _{a,b,\beta }^{\ }\Vert z\Vert _{a,b,\beta }^{\ }\right)
(b-a)^{2\beta } \\
&&+k\left( \left\| x\right\| _{a,b,\beta }\left\| y\right\| _{a,b,\beta }^{\
}\Vert z\Vert _{a,b,\beta }^{\ }+\left\| x\right\| _{a,b,\beta }\left\|
y\otimes z\right\| _{a,b,\beta }+\left\| z\right\| _{a,b,\beta }\left\|
x\otimes y\right\| _{a,b,\beta }\right)  \\
&&\times \left( \Vert f^{\prime }\Vert _{\infty }+\Vert f^{\prime }\Vert _{{%
\lambda }}\Vert x\Vert _{a,b,\beta }^{{\lambda }}(b-a)^{\lambda \beta
}\right) ](b-a)^{3{\beta -2\varepsilon }},
\end{eqnarray*}%
which implies the desired result.
\end{proof}

 \setcounter{equation}{0}
\section{Differential Equations Driven by Rough Paths}
Let $y:[0,1]\rightarrow \mathbb{R}^{d}$ be a $\beta $-H\"{o}lder continuous
function. Suppose that \ $(y^{i},y^{j},y^{i}\otimes y^{j})$ is a $\beta $-H%
\"{o}lder continuous multiplicative function, for each $i,j=1,\ldots ,d$. We
aim to solve the differential equation%
\begin{equation}
x_{t}=x_{0}+\int_{0}^{t}f(x_{r})dy_{r},  \label{3.1}
\end{equation}%
where $f=\mathbb{R}^{m}\rightarrow \mathbb{R}^{md}$.

Formula (\ref{2.9}) and Definition \ref{d.3.3} allow us to transform this
equation into the following system of integral equations:

\begin{eqnarray}
x_{t} &=&x_{0}+(-1)^{\alpha }\int_{0}^{t}\widehat{D}_{0+}^{\alpha }f\left(
x\right) (s)D_{t-}^{1-\alpha }y_{t-}(s)ds  \label{3.2} \\
&&+\int_{0}^{t}D_{0+}^{2\alpha -1}f^{\prime }\left( x\right)
(s)\int_{s}^{t}\int_{s}^{\eta }K_{0,s}(\xi ,\eta )\Gamma ^{\alpha -{%
\varepsilon }}\left( x\otimes y\right) _{\xi ,\eta }d\xi d\eta ds,  \notag
\end{eqnarray}%
\begin{eqnarray}
\left( x\otimes y\right) _{s,t} &=&(-1)^{\alpha }\int_{s}^{t}\widehat{D}%
_{s+}^{\alpha }f\left( x\right) (r)D_{t-}^{1-\alpha }\left( y\otimes
y\right) _{\cdot ,t-}(r)dr  \notag \\
&&+\int_{s}^{t}D_{s+}^{2\alpha -1}f^{\prime }\left( x\right) (r)  \label{3.3}
\\
&&\times \int_{r}^{t}\int_{r}^{\eta }K_{s,r}(\xi ,\eta )\Gamma ^{\alpha
-\varepsilon }\left( x\otimes \left( y\otimes y\right) _{\cdot ,t}\right)
_{\xi ,\eta }d\xi d\eta dr.  \notag
\end{eqnarray}

\begin{theorem}
\label{th4} Let $y:[0,1]\rightarrow \mathbb{R}^{d}$ be a $\beta $-H\"{o}lder
continuous function. Suppose that \ $(y^{i},y^{j},y^{i}\otimes y^{j})$ is a
real valued $\beta $-H\"{o}lder continuous multiplicative function, for each 
$i,j=1,\ldots ,d$. Let $f:\mathbb{R}^{m}\rightarrow \mathbb{R}^{md}$ be a
continuously differentiable function such that $f^{\prime }$ is $\lambda $-H%
\"{o}lder continuous, where $\lambda >\frac{1}{\beta }-2$, and $f$ and $%
f^{\prime }$ are bounded.   Set%
\begin{equation*}
\rho _{f}:=\Vert f\Vert _{\infty }+\Vert f^{\prime }\Vert _{\infty }+\Vert
f^{\prime }\Vert _{{\lambda }}.
\end{equation*}%
Then there is a solution $\ $to Equations (\ref{3.2})- (\ref{3.3}), such
that $(x,y,x\otimes y)$ is a $\beta $-H\"{o}lder continuous multiplicative
functional. Moreover,  for any $\gamma >\frac 1\beta$ the function $x$ satisfies the estimate%
\begin{equation}
\sup_{0\leq t\leq T}|x_{t}|\leq |x_{0}|+\ T\left\{ 2k\rho _{f}\left[ \Vert
y\Vert _{{\beta }}+\frac{\Vert y\otimes y\Vert _{\beta }}{\Vert y\Vert
_{\beta }}\right] \vee 1\right\} ^\gamma\,,
\label{3.4}
\end{equation}%
where $k$ is a universal constant depending only on  $\beta $ and $\gamma$.
\end{theorem}

\begin{proof}
To simplify the proof we will assume $d=m=1$. The proof will be done in
several steps.

\textbf{Step 1.}   Fix $\alpha>0$ and $\varepsilon>0$ such that
$1-\beta <\alpha< 2\beta$, $\alpha <\frac{\lambda\beta +1}2$,
$\varepsilon< \alpha+\beta-1$, $\varepsilon <\frac \beta 2$, and
${(1-2\varepsilon )/({\beta -2\varepsilon )}}<\gamma$. 

We will write the Equations (\ref{3.2}) and (\ref{3.3}) in
the compact form%
\begin{eqnarray*}
x &=&\Phi _{1}(x,y,x\otimes y), \\
x\otimes y &=&\Phi _{2}(x,y,y\otimes y,x\otimes y).
\end{eqnarray*}%
Consider the mapping $J:(x,x\otimes y)\rightarrow (J_{1}x,J_{2}\left(
x\otimes y\right) )$ defined by%
\begin{eqnarray*}
J_{1}x &=&\Phi _{1}(x,y,x\otimes y), \\
J_{2}\left( x\otimes y\right) &=&\Phi _{2}(x,y,y\otimes y,x\otimes y).
\end{eqnarray*}%
We need some a priori estimates of the H\"{o}lder norms of $J_{1}x$ and $%
J_{2}\left( x\otimes y\right) $ in terms of the H\"{o}lder norms of $x$ and $%
x\otimes y$. From (\ref{2.15}) it follows that 
\begin{eqnarray}
\left\| J_{1}x\right\| _{s,t,\beta } &\leq &k[\left\| f\right\| _{\infty
}\left\| y\right\| _{s,t,\beta }+\left( \left\| x\otimes y\right\|
_{s,t,\beta }+\left\| x\right\| _{s,t,\beta }\left\| y\right\| _{s,t,\beta
}\right)  \notag \\
&&\qquad \times \left( \left\| f^{\prime }\right\| _{\infty }+\left\|
f^{\prime }\right\| _{\lambda }\left\| x\right\| _{s,t,\beta }^{\lambda
}(t-s)^{{\lambda }{\beta }}\right) (t-s)^{\beta -2\varepsilon }].
\label{3.5}
\end{eqnarray}%
On the other hand, Proposition \ref{p.3.1} implies that 
\begin{equation}
\left\| J_{2}\left( x\otimes y\right) \right\| _{s,t,\beta }\leq \ k\left[
A_{s,t}+B_{s,t}\left\| x\otimes y\right\| _{s,t,{\beta }}\,(t-s)^{\beta
-2\varepsilon }\right] ,  \label{3.6}
\end{equation}%
where $A_{s,t}$ and $B_{s,t}$ are defined by (\ref{2.29}) and (\ref{2.30}),
respectively.

\textbf{Step 2. \ }Set 
\begin{equation*}
\alpha (y):=\left( \ 2k\rho _{f}\left[ \Vert y\Vert _{{\beta }}+\frac{\Vert
y\otimes y\Vert _{\beta }}{\Vert y\Vert _{\beta }}\right]  \vee 1\right) ^{1/(\beta
-2\varepsilon )},
\end{equation*}%
where $k$ is the constant appearing in formulas (\ref{3.5}) and (\ref{3.6}).
Suppose that 
\begin{equation}
0<t-s\leq \frac{1}{\alpha (y)}.  \label{3.7}
\end{equation}%
Then, the inequalities 
\begin{eqnarray}
\Vert x\Vert _{s,t,\beta } &\leq &2k\rho _{f}\Vert y\Vert _{\beta }
\label{3.8} \\
(t-s)^{\beta -2\varepsilon }\left\| x\otimes y\right\| _{s,t,{\beta }} &\leq
&\ \Vert y\Vert _{\beta }  \label{3.9}
\end{eqnarray}%
imply that 
\begin{eqnarray}
\Vert J_{1}x\Vert _{s,t,\beta } &\leq &2k\rho _{f}\Vert y\Vert _{\beta }
\label{3.10} \\
(t-s)^{\beta -2\varepsilon }\left\| J_{2}(x\otimes y)\right\| _{s,t,{\beta }%
} &\leq &\ \Vert y\Vert _{\beta }  \label{3.11}
\end{eqnarray}%
In fact, from the definition of $\alpha (y)$ and (\ref{3.8}) we deduce%
\begin{equation}
(t-s)^{\beta -2\varepsilon }\Vert x\Vert _{s,t,{\beta }}\leq 1.  \label{3.12}
\end{equation}

By the definition of $B_{s,t}$ and $A_{s,t}$ we have 
\begin{equation}
B_{s,t}\leq (\Vert f^{\prime }\Vert _{\infty }+\Vert f^{\prime }\Vert _{{%
\lambda }})\Vert y\Vert _{\beta }\,\leq \rho _{f}\Vert y\Vert _{\beta }\leq 
\frac{(t-s)^{-(\beta -2\varepsilon )}}{2k},  \label{3.13}
\end{equation}%
and%
\begin{eqnarray}
A_{s,t} &\leq &\left( \Vert f\Vert _{\infty }+\Vert f^{\prime }\Vert
_{\infty }+\Vert f^{\prime }\Vert _{{\lambda }}\right) \left( \left\|
y\otimes y\right\| _{\beta }+\left\| y\right\| _{\beta }^{2}\right) \, 
\notag \\
&\leq &\frac{(t-s)^{-(\beta -2\varepsilon )}}{2k}\Vert y\Vert _{\beta }.
\label{3.14}
\end{eqnarray}%
Therefore, substituting (\ref{3.13}) and (\ref{3.14}) into (\ref{3.6}) we
obtain (\ref{3.11}). Finally, from (\ref{3.5}) we get (\ref{3.10}).

\textbf{Step 3.} We can now proceed with the proof of the existence. Let $N$
be a natural number such that $\frac{T}{N}=\delta \leq \frac{1}{\alpha (y)}$. 
We partition the interval $[0,T]$ in $N$ subintervals of the same
length and set $t_{i}=\frac{iT}{N}$, $i=0,1,\ldots ,N-1$. We will make use
of the notation $\left\| x\right\| _{i}=\Vert x\Vert _{t_{i-1},t_{i},\beta }$%
, and $\left\| x\otimes y\right\| _{i}=\Vert x\otimes y\Vert
_{t_{i-1},t_{i},\beta }$, for $i=1,\ldots ,N-1$. \ From Step 2 we know that
if that $x$ and $x\otimes y$ satisfy 
\begin{eqnarray*}
\Vert x\Vert _{i} &\leq &\ 2k\rho _{f}\Vert y\Vert _{\beta } \\
\left\| x\otimes y\right\| _{i} &\leq &\ \Vert y\Vert _{\beta }\delta
^{-\left( \beta -2\varepsilon \right) },
\end{eqnarray*}%
for any $i=1,\ldots ,N-1$, then the same inequalities hold for $Jx$ and $%
Jx\otimes y$, that is%
\begin{eqnarray*}
\Vert J_{1}x\Vert _{i} &\leq &\ 2k\rho _{f}\Vert y\Vert _{\beta } \\
\left\| J_{2}(x\otimes y)\right\| _{i} &\leq &\ \ \Vert y\Vert _{\beta
}\delta ^{-\left( \beta -2\varepsilon \right) }.
\end{eqnarray*}%
Consequently, there is a constant $C_{1}$ such that 
\begin{equation*}
\Vert J_{1}^{n}x\Vert _{\beta }+\Vert J_{2}^{n}(x\otimes y)\Vert _{\beta
}\leq C_{1}\,.
\end{equation*}%
This implies that the sequence of functions $J_{1}^{n}x$ is equicontinuous
and bounded in $C^{\beta }$. Therefore, there exists a subsequence which
converges in the $\beta ^{\prime }$-H\"{o}lder norm if $\beta ^{\prime
}<\beta $. In the same way, there is a subsequence of $J_{2}^{n}(x\otimes y)$
which converges in the $\beta ^{\prime }$-H\"{o}lder norm. The limit $%
(x,x\otimes y)$ defines a $\beta $-H\"{o}lder continuous multiplicative
functional $(x,y,x\otimes y)$. Using the continuity of the solution in this
norm it is not difficult to show that the limit is a solution. This implies
the existence of a solution, which satisfies (\ref{3.8}) and (\ref{3.9}).

\textbf{Step 4.} \ Let us now prove the estimate (\ref{3.4}). By step 2, the
solution we have constructed satisfies the estimates (\ref{3.8}) and (\ref%
{3.9}) if (\ref{3.7}) holds. Then it follows that for any $r\in \lbrack s,t]$%
\begin{equation*}
\sup_{r\in \lbrack s,t]}|x_{r}|\leq |x_{s}|+(t-s)^{\beta }\Vert x\Vert
_{s,t,\beta }\leq |x_{s}|+(t-s)^{2\varepsilon }.
\end{equation*}%
Since the interval $[0,T]$ can be divided into $[T/\tau ]$ intervals of
length $\tau =\frac{1}{\alpha (y)}$, the inequality (\ref{3.4}) follows.
\end{proof}

\begin{theorem}
\label{th5} Let $y:[0,1]\rightarrow \mathbb{R}^{d}$ be a $\beta $-H\"{o}lder
continuous function. Suppose that \ $(y^{i},y^{j},y^{i}\otimes y^{j})$ is a
real valued $\beta $-H\"{o}lder continuous multiplicative function, for each 
$i,j=1,\ldots ,d$.   Let $f:\mathbb{R}^{m}\rightarrow \mathbb{R}^{md}$ be a
twice continuously differentiable function such  that $f^{\prime \prime }$ is $%
\lambda $-H\"{o}lder continuous, where $\lambda >\frac{1}{\beta }-2$, and $f$%
, $f^{\prime }$ and $f^{\prime \prime }$ are bounded.   Then there is a unique solution to Equations
(\ref{3.2})- (\ref{3.3}) such that $(x,y,x\otimes y)$ is a $\beta $-H\"{o}%
lder continuous multiplicative functional.

Moreover, if $\tilde{x}$ satisfies $\tilde{x}_{t}=\tilde{x}%
_{0}+\int_{0}^{t}f(\widetilde{x}_{r})d\tilde{y}_{r}$\ and $\tilde{y}$ verifies the same
hypotheses as $y$, then 
\begin{equation}
\sup_{0\leq t\leq T}|x_{t}-\tilde{x}_{t}|\leq 
 C\left\{ |x_{0}-\tilde{x}_{0}|+\Vert y-\tilde{y}%
\Vert _{{\beta }}+\Vert y\otimes (y-\tilde{y})\Vert _{{\beta }}\right\} \,,
\label{3.15}
\end{equation}%
where $C$ depends on $\|y\|_\beta$, $\|y\otimes y\|_\beta$,
$\beta$,   $\lambda$, and $\hat{\rho}_f$, and where  
\begin{equation*}
\hat{\rho}_{f}=\Vert f\Vert _{\infty }+\left\| f\right\| _{\lambda }+\Vert
f^{\prime }\Vert _{\infty }+\left\| f^{\prime }\right\| _{\lambda }+\Vert
f^{\prime \prime }\Vert _{\infty }+\Vert f^{\prime \prime }\Vert _{{\lambda }%
}\,.
\end{equation*}
\end{theorem}

\begin{proof}
To simplify the proof we will assume $d=m=1$. Notice that uniqueness follows
from the estimate (\ref{3.15}). So it suffices to show this inequality. We
fix \ $s<t$ such that   $t-s\leq \frac 1 {\beta(y)}  $, where $\beta$ is defined as follows
\begin{equation}
\beta(y)= \left( \ 2k\widehat{\rho }_{f}\left[ \Vert y\Vert _{{\beta }}+
\|y\|^2_\beta+ \Vert y\otimes y\Vert _{\beta }
+
\frac{\Vert y\otimes y\Vert _{\beta }} {\Vert y\Vert _{\beta }}\right] \vee 1
\right) ^{1/(\beta \lambda )}\,.   \label{beta}
\end{equation}%
The constant $k$ appearing in the definition of $\beta$ will be chosen later.
We choose $\alpha>0$ and $\varepsilon>0$ such that
$1-\beta <\alpha< 2\beta$, $\alpha <\frac{\lambda\beta +1}2$, and
$\varepsilon< \alpha+\beta-1$, $\varepsilon <\frac \beta 2$.
 We also assume that the solutions $x$ and $\widetilde{x}$ satisfy the following inequalities:  
\begin{equation}
(t-s) ^{\beta -2\varepsilon }\Vert x\Vert _{{\beta }}\leq 1,  \label{3.16}
\end{equation}%
\begin{equation}
(t-s) ^{\beta -2\varepsilon }\Vert \widetilde{x}\Vert _{{\beta }}\leq 1,
\label{3.17}
\end{equation}%
\begin{equation}
(t-s) ^{\beta -2\varepsilon }\left\| x\otimes y\right\| _{{\beta }}\leq
\Vert y\Vert _{\beta }.  \label{3.18}
\end{equation}%

Our first purpose is to estimate the H\"{o}lder norm $\left\| x-\widetilde{x}%
\right\| _{s,t,\beta }$. We can write%
\begin{eqnarray*}
\left\| x-\widetilde{x}\right\| _{s,t,\beta } &\leq &\left\| \int \left[
f(x_{s})-f(\tilde{x}_{s})\right] dy_{s}\right\| _{s,t,\beta }+\left\| \int f(%
\tilde{x}_{s})d(y_{s}-\tilde{y}_{s})\right\| _{s,t,\beta } \\
&=&I_{1,s,t}+I_{2,s,t}.
\end{eqnarray*}%
The term $I_{1,s,t}$ can be estimated using (\ref{2.16}) and we obtain 
\begin{equation}
I_{1,s,t}\leq k[H_{1}\Vert x-\tilde{x}\Vert _{s,t,\infty }+H_{2}\Vert x-%
\tilde{x}\Vert _{s,t,\beta }+H_{3}\Vert (x-\tilde{x})\otimes y\Vert _{s,t,{%
\beta }}],  \label{3.19}
\end{equation}%
where%
\begin{eqnarray*}
H_{1} &=&\left\| y\right\| _{s,t,\beta }\left( \Vert f^{\prime }\Vert
_{\infty }+\left\| f^{\prime \prime }\right\| _{\lambda }\Vert \widetilde{x}%
\Vert _{s,t,{\beta }}\left( \Vert x\Vert _{s,t,{\beta }}^{\lambda }+\Vert 
\tilde{x}\Vert _{s,t,{\beta }}^{\lambda }\right) \right) (t-s)^{\beta
(1+\lambda )} \\
&&+\left( \Vert f^{\prime \prime }\Vert _{\infty }+\Vert f^{\prime \prime
}\Vert _{\lambda }\left( \Vert x\Vert _{s,t,\beta }^{\lambda }+\Vert 
\widetilde{x}\Vert _{s,t,\beta }^{\lambda }\right) (t-s)^{\beta \lambda
}\right)  \\
&&\times \left( \Vert x\otimes y\Vert _{s,t,{\beta }}+\Vert x\Vert _{s,t,{%
\beta }}\left\| y\right\| _{s,t,\beta }\ \right) (t-s)^{\beta -2\varepsilon
},
\end{eqnarray*}%
\begin{eqnarray*}
H_{2} &=&\ \left\| f^{\prime \prime }\right\| _{\infty }\left\| y\right\|
_{s,t,\beta }\left( \Vert x\Vert _{s,t,{\beta }}+\Vert \tilde{x}\Vert _{s,t,{%
\beta }}\right) (t-s)^{\beta (1+\lambda )}\  \\
&&+\left\| f^{\prime \prime }\right\| _{\infty }\left( \Vert x\otimes y\Vert
_{s,t,\beta }+\Vert x\Vert _{s,t,\beta }\left\| y\right\| _{s,t,\beta }\
\right) (t-s)^{2\beta -2\varepsilon } \\
&&+\ \left( \Vert f\Vert _{\infty }+\Vert f\Vert _{\lambda }\Vert \tilde{x}%
\Vert _{s,t,\beta }^{\lambda }(t-s)^{\lambda \beta }\right) \left\|
y\right\| _{s,t,\beta }(t-s)^{\beta -2\varepsilon },
\end{eqnarray*}%
and%
\begin{equation*}
H_{3}=\left( \Vert f\Vert _{\infty }+\Vert f\Vert _{\lambda }\Vert \tilde{x}%
\Vert _{s,t,\beta }^{\lambda }(t-s)^{\lambda \beta }\right) (t-s)^{\beta
-2\varepsilon }.
\end{equation*}%
Then, using the inequalities (\ref{3.16}), (\ref{3.17}), and (\ref{3.18}) we
get the following estimates%
\begin{eqnarray}
H_{1} &\leq &\left\| y\right\| _{\beta }\left( \Vert f^{\prime }\Vert
_{\infty }+2\Vert f^{\prime \prime }\Vert _{\infty }+6\left\| f^{\prime
\prime }\right\| _{\lambda }\right), \   \label{3.20} \\
H_{2} &\leq &\ \left\| y\right\| _{\beta }\left( 3\left\| f^{\prime \prime
}\right\| _{\infty }+\ \Vert f\Vert _{\infty }+\Vert f\Vert _{\lambda
}\right) (t-s)^{\beta \lambda },  \label{3.21} \\
H_{3} &\leq &\left( \Vert f\Vert _{\infty }+\Vert f\Vert _{\lambda }\right) .
\label{3.22}
\end{eqnarray}%
It remains to handle the term $\Vert (x-\tilde{x})\otimes y\Vert _{s,t,{%
\beta }}$ in (\ref{3.19}). To get estimates for this term we apply   again
the inequality (\ref{2.16}) and we have 
\begin{eqnarray}
\left| ((x-\tilde{x})\otimes y)_{s,t}\right|  &=&\left| \int_{s}^{t}\left[
f(x_{r})-f(\tilde{x}_{r})\right] d_{r}(y\otimes y)_{r,t}\right|   \notag \\
&\leq &k(t-s)^{\beta }\left[ \widetilde{H}_{1}\Vert x-\tilde{x}\Vert
_{s,t,\infty }+\widetilde{H}_{2}\Vert x-\tilde{x}\Vert _{s,t,\beta }\right. 
\notag \\
&&\left. +\widetilde{H}_{3}\Vert (x-\tilde{x})\otimes (y\otimes y)_{\cdot
,t}\Vert _{s,t,{\beta }}\right] ,  \label{3.23}
\end{eqnarray}%
where%
\begin{eqnarray*}
\widetilde{H}_{1} &=&\left\| (y\otimes y)_{\cdot ,t}\right\| _{s,t,\beta
}\left( \Vert f^{\prime }\Vert _{\infty }+\left\| f^{\prime \prime }\right\|
_{\lambda }\Vert \widetilde{x}\Vert _{s,t,{\beta }}\left( \Vert x\Vert _{s,t,%
{\beta }}^{\lambda }+\Vert \tilde{x}\Vert _{s,t,{\beta }}^{\lambda }\right)
\right) (t-s)^{\beta (1+\lambda )} \\
&&+\left( \Vert f^{\prime \prime }\Vert _{\infty }+\Vert f^{\prime \prime
}\Vert _{\lambda }\left( \Vert x\Vert _{s,t,\beta }^{\lambda }+\Vert 
\widetilde{x}\Vert _{s,t,\beta }^{\lambda }\right) (t-s)^{\beta \lambda
}\right)  \\
&&\times \left( \Vert x\otimes (y\otimes y)_{\cdot ,t}\Vert _{s,t,{\beta }%
}+\Vert x\Vert _{s,t,{\beta }}\left\| (y\otimes y)_{\cdot ,t}\right\|
_{s,t,\beta }\ \right) (t-s)^{\beta -2\varepsilon },
\end{eqnarray*}%
\begin{eqnarray*}
\widetilde{H}_{2} &=&\ \left\| f^{\prime \prime }\right\| _{\infty }\left\|
(y\otimes y)_{\cdot ,t}\right\| _{s,t,\beta }\left( \Vert x\Vert _{s,t,{%
\beta }}+\Vert \tilde{x}\Vert _{s,t,{\beta }}\right) (t-s)^{\beta (1+\lambda
)}\  \\
&&+\left\| f^{\prime \prime }\right\| _{\infty }\left( \Vert x\otimes
(y\otimes y)_{\cdot ,t}\Vert _{s,t,\beta }+\Vert x\Vert _{s,t,\beta }\left\|
(y\otimes y)_{\cdot ,t}\right\| _{s,t,\beta }\ \right) (t-s)^{2\beta
-2\varepsilon } \\
&&+\ \left( \Vert f\Vert _{\infty }+\Vert f\Vert _{\lambda }\Vert \tilde{x}%
\Vert _{s,t,\beta }^{\lambda }(t-s)^{\lambda \beta }\right) \left\|
(y\otimes y)_{\cdot ,t}\right\| _{s,t,\beta }(t-s)^{\beta -2\varepsilon },
\end{eqnarray*}%
and%
\begin{equation*}
\widetilde{H}_{3}=H_{3}.
\end{equation*}%
Using (\ref{2.31}), (\ref{2.32}), (\ref{3.16}), (\ref{3.17}) and (\ref{3.18}%
) we get the following estimates%
\begin{eqnarray}
\lefteqn{\widetilde{H}_{1} =\left( \left\| y\otimes y\right\| _{s,t,\beta }+\
\left\| y\right\| _{s,t,\beta }^{2}\right)  } \notag \\
&& \times\left( \Vert f^{\prime }\Vert
_{\infty }+\left\| f^{\prime \prime }\right\| _{\lambda }\Vert \widetilde{x}%
\Vert _{s,t,{\beta }}\left( \Vert x\Vert _{s,t,{\beta }}^{\lambda }+\Vert 
\tilde{x}\Vert _{s,t,{\beta }}^{\lambda }\right) \right) (t-s)^{\beta
(2+\lambda )}  \notag \\
&&+k\left( \Vert f^{\prime \prime }\Vert _{\infty }+\Vert f^{\prime \prime
}\Vert _{\lambda }\left( \Vert x\Vert _{s,t,\beta }^{\lambda }+\Vert 
\widetilde{x}\Vert _{s,t,\beta }^{\lambda }\right) (t-s)^{\beta \lambda
}\right)   \notag \\
&&\times \left( \left\| x\right\| _{s,t\beta }\left\| y\right\| _{s,t,\beta
}^{2}+\left\| x\right\| _{s,t,\beta }\left\| y\otimes y\right\| _{s,t,\beta
}+\left\| y\right\| _{s,t,\beta }\left\| x\otimes y\right\| _{s,t,\beta
}\right) \ (t-s)^{2\beta -2\varepsilon }  \notag \\
&\leq &k\left( \left\| y\otimes y\right\| _{\beta }+\ \left\| y\right\|
_{\beta }^{2}\right) \left( \Vert f^{\prime }\Vert _{\infty }+\Vert
f^{\prime \prime }\Vert _{\infty }+\left\| f^{\prime \prime }\right\|
_{\lambda }\right) (t-s)^{\beta },  \label{3.24}
\end{eqnarray}%
and%
\begin{eqnarray}
\lefteqn{\widetilde{H}_{2} =\ \left\| f^{\prime \prime }\right\| _{\infty }\left(
\left\| y\otimes y\right\| _{s,t,\beta }+\ \left\| y\right\| _{s,t,\beta
}^{2}\right) \left( \Vert x\Vert _{s,t,{\beta }}+\Vert \tilde{x}\Vert _{s,t,{%
\beta }}\right) (t-s)^{\beta (2+\lambda )}\  } \notag \\
&&+k\left\| f^{\prime \prime }\right\| _{\infty }\Big( \left\| y\right\|
_{s,t,\beta }\left\| x\otimes y\right\| _{s,t,\beta }   \\
&&+\Vert x\Vert
_{s,t,\beta }\left\| y\otimes y\right\| _{s,t,\beta }+\Vert x\Vert
_{s,t,\beta }\ \left\| y\right\| _{s,t,\beta }^{2}\ \Big) (t-s)^{3\beta
-2\varepsilon }  \notag \\
&&+k\left( \Vert f\Vert _{\infty }+\Vert f\Vert _{\lambda }\Vert \tilde{x}%
\Vert _{s,t,\beta }^{\lambda }(t-s)^{\lambda \beta }\right) \left( \left\|
y\otimes y\right\| _{s,t,\beta }+\ \left\| y\right\| _{s,t,\beta
}^{2}\right) (t-s)^{2\beta -2\varepsilon }  \notag \\
&\leq &\ k\left( \Vert f\Vert _{\infty }+\Vert f\Vert _{\lambda }+\left\|
f^{\prime \prime }\right\| _{\infty }\right) \left( \left\| y\otimes
y\right\| _{\beta }+\ \left\| y\right\| _{\beta }^{2}\right) (t-s)^{\beta }.
\label{3.25}
\end{eqnarray}%
On the other hand, from (\ref{2.32}) we get 
\begin{eqnarray}
\lefteqn{\Vert (x-\tilde{x})\otimes (y\otimes y)_{\cdot ,t}\Vert _{s,t,{\beta }}
\leq k\Big( \left\| x-\tilde{x}\right\| _{s,t\beta }\left\| y\right\|
_{\beta }^{2}  }\notag  \\ \label{3.26} 
&&+\left\| x-\tilde{x}\right\| _{s,t,\beta }\left\| y\otimes
y\right\| _{\beta }+\left\| y\right\| _{\beta }\left\| \left( x-%
\tilde{x}\right) \otimes y\right\| _{s,t,\beta }\Big) (t-s)^{{\beta }}. 
\end{eqnarray}%
Thus, substituting (\ref{3.24}), (\ref{3.25}), (\ref{3.22}) and (\ref{3.26})
into (\ref{3.23}) yields%
\begin{eqnarray*}
\lefteqn{\left\| (x-\tilde{x})\otimes y\right\| _{s,t,\beta } \leq k(t-s)^{\beta
}\Big[\left( \left\| y\otimes y\right\| _{\beta }+\ \left\| y\right\| _{\beta
}^{2}\right) }  \\
&&\times \left( \Vert f^{\prime }\Vert _{\infty }+\Vert f^{\prime \prime
}\Vert _{\infty }+\left\| f^{\prime \prime }\right\| _{\lambda }\right)
\Vert x-\tilde{x}\Vert _{s,t,\infty } \\
&&+\left( \Vert f\Vert _{\infty }+\Vert f\Vert _{\lambda }+\left\| f^{\prime
\prime }\right\| _{\infty }\right) \left( \left\| y\otimes y\right\| _{\beta
}+\ \left\| y\right\| _{\beta }^{2}\right) \Vert x-\tilde{x}\Vert
_{s,t,\beta } \\
&&+\left( \Vert f\Vert _{\infty }+\Vert f\Vert _{\lambda }\right)  \\
&&\times \left( \left\| x-\tilde{x}\right\| _{s,t\beta }\left\| y\right\|
_{\beta }^{2}+\left\| x-\tilde{x}\right\| _{s,t,\beta }\left\| y\otimes
y\right\| _{\beta }+\left\| y\right\| _{\beta }\left\| \left( x-\tilde{x}%
\right) \otimes y\right\| _{s,t,\beta }\right) \Big] \\
&\leq &k(t-s)^{\beta }\hat{\rho}_{f}\left( \left\| y\otimes y\right\|
_{\beta }+\ \left\| y\right\| _{\beta }^{2}\right) \left( \Vert x-\tilde{x}%
\Vert _{s,t,\infty }+\Vert x-\tilde{x}\Vert _{s,t,\beta }\right)  \\
&&+k(t-s)^{\beta }\left( \Vert f\Vert _{\infty }+\Vert f\Vert _{\lambda
}\right) \left\| y\right\| _{\beta }\left\| \left( x-\tilde{x}\right)
\otimes y\right\| _{s,t,\beta }.
\end{eqnarray*}%
The condition $t-s\le 1/\beta(y)$, if the constant in $\beta(y)$ is chosen in an appropriate
way, implies that
$$
k(t-s)^{\ \beta }\left( \Vert
f\Vert _{\infty }+\Vert f\Vert _{\lambda }\right) \left\| y\right\| _{\beta
}\leq \frac{1}{2}.
$$
 Hence,  
\begin{equation}
\left\| (x-\tilde{x})\otimes y\right\| _{s,t,\beta }\leq k(t-s)^{\ \beta }%
\hat{\rho}_{f}\left( \left\| y\otimes y\right\| _{\beta }+\ \left\|
y\right\| _{\beta }^{2}\right) \left( \Vert x-\tilde{x}\Vert _{s,t,\infty
}+\Vert x-\tilde{x}\Vert _{s,t,\beta }\right) .  \label{3.27}
\end{equation}%
Substituting (\ref{3.27}), (\ref{3.20}), (\ref{3.21}) and (\ref{3.22}) into (%
\ref{3.19}) yields%
\begin{eqnarray*}
I_{1,s,t} &\leq &k\hat{\rho}_{f}[\left\| y\right\| _{\beta }\Vert x-\tilde{x}%
\Vert _{s,t,\infty }+\left\| y\right\| _{\beta }\Vert x-\tilde{x}\Vert
_{s,t,\beta }(t-s)^{\beta \lambda } \\
&&+\left( \left\| y\otimes y\right\| _{\beta }+\ \left\| y\right\| _{\beta
}^{2}\right) \left( \Vert x-\tilde{x}\Vert _{s,t,\infty }+\Vert x-\tilde{x}%
\Vert _{s,t,\beta }\right) (t-s)^{\beta }].
\end{eqnarray*}%
Again, condition $t-s\le 1/\beta(y)$, if the constant in $\beta(y)$ is chosen in an appropriate
way, implies that
\begin{equation}
I_{1,s,t}\leq k\hat{\rho}_{f}\left\| y\right\| _{\beta }\ \Vert x-\tilde{x}%
\Vert _{s,t,\infty }+\frac{1}{2}\Vert x-\tilde{x}\Vert _{s,t,\beta }\ .
\label{3.28}
\end{equation}

For the term $I_{2,s,t}$ we have the following estimates, using (\ref{2.25})%
\begin{eqnarray}
I_{2,s,t} &\leq &k\left\| f\right\| _{\infty }\left\| y-\tilde{y}\right\|
_{s,t,\beta }+k\left( \left\| x\otimes (y-\tilde{y})\right\| _{s,t,\beta
}+\left\| x\right\| _{s,t,\beta }\left\| y-\tilde{y}\right\| _{s,t,\beta
}\right)   \notag \\
&&\times \left( \left\| f^{\prime }\right\| _{\infty }+\left\| f^{\prime
}\right\| _{\lambda }\left\| x\right\| _{s,t,\beta }^{\lambda
}(t-s)^{\lambda \beta }\right) (t-s)^{\beta -2\varepsilon }  \notag \\
&\leq &k\rho _{f}\left\| y-\tilde{y}\right\| _{\beta }+k\left\| x\otimes (y-%
\tilde{y})\right\| _{s,t,\beta }\left( \left\| f^{\prime }\right\| _{\infty
}+\left\| f^{\prime }\right\| _{\lambda }\right) (t-s)^{\beta -2\varepsilon
}.  \label{3.29}
\end{eqnarray}%
In order to estimate $\left\| x\otimes (y-\tilde{y})\right\| _{s,t,\beta }$
we make use of Proposition \ \ref{p.3.1} and we obtain%
\begin{eqnarray}
\left\| x\otimes (y-\tilde{y})\right\| _{s,t,\beta } &=&\sup_{s\leq \xi \leq
\eta \leq t}\frac{1}{(\eta -\xi )^{2\beta }}\left| \int_{\xi }^{\eta
}f(x_{r})d_{r}(y\otimes (y-\tilde{y}))_{r,\eta }\right|   \notag \\
&\leq &\ k\left[ A_{s,t}+B_{s,t}\left\| x\otimes y\right\| _{s,t,\beta
}\,(t-s)^{\beta -2\varepsilon }\right] ,  \label{3.31}
\end{eqnarray}%
where 
\begin{eqnarray}
A_{s,t} &=&\left( \Vert y\otimes (y-\tilde{y})\Vert _{s,t,\beta }+\Vert
y\Vert _{s,t,\beta }^{\ }\Vert (y-\tilde{y})\Vert _{s,t,\beta }^{\ }\right) 
\notag \\
&&\times \left[ \left\| f\right\| _{\infty }+\left( \left\| x\right\|
_{s,t,\beta }\left\| f^{\prime }\right\| _{\infty }+\left\| f^{\prime
}\right\| _{\lambda }\left\| x\right\| _{s,t,\beta }^{1+\lambda
}(t-s)^{\lambda \beta }\right) (t-s)^{\beta -2\varepsilon }\right] \,  \notag
\\
&\leq &\rho _{f}\left( \Vert y\otimes (y-\tilde{y})\Vert _{s,t,\beta }+\Vert
y\Vert _{s,t,\beta }^{\ }\Vert (y-\tilde{y})\Vert _{s,t,\beta }^{\ }\right) ,
\label{3.33}
\end{eqnarray}%
and 
\begin{eqnarray}
B_{s,t} &=&\ \left\| (y-\tilde{y})\right\| _{s,t,\beta }\left( \left\|
f^{\prime }\right\| _{\infty }+\left\| f^{\prime }\right\| _{\lambda
}\left\| x\right\| _{s,t,\beta }^{\lambda }(t-s)^{\lambda \beta }\right) \, 
\notag \\
&\leq &\left\| (y-\tilde{y})\right\| _{s,t,\beta }\left( \left\| f^{\prime
}\right\| _{\infty }+\left\| f^{\prime }\right\| _{\lambda }\right) .
\label{3.34}
\end{eqnarray}%
Substituting (\ref{3.33}) and (\ref{3.34}) into (\ref{3.31}) yields 
\begin{equation}
\left\| x\otimes (y-\tilde{y})\right\| _{s,t,\beta }\leq \ k\rho _{f}\left[
\Vert y\otimes (y-\tilde{y})\Vert _{\beta }+\Vert y\Vert _{\beta }^{\ }\Vert
(y-\tilde{y})\Vert _{\beta }^{\ }\right] .  \label{3.35}
\end{equation}%
Finally, from (\ref{3.35}) and (\ref{3.39}) we obtain%
\begin{eqnarray}
I_{2,s,t} &\leq &k\rho _{f}\left\| y-\tilde{y}\right\| _{\beta }+k\rho
_{f}^{2}\left[ \Vert y\otimes (y-\tilde{y})\Vert _{\beta }+\Vert y\Vert
_{s,t,\beta }^{\ }\Vert (y-\tilde{y})\Vert _{\beta }^{\ }\right]
(t-s)^{\beta -2\varepsilon }  \notag \\
&\leq &k\rho _{f}\left\| y-\tilde{y}\right\| _{\beta }+k\rho _{f}^{2}\Vert
y\otimes (y-\tilde{y})\Vert _{\beta }(t-s)^{\beta -2\varepsilon }.
\label{3.37}
\end{eqnarray}%
Now from (\ref{3.28}) and (\ref{3.37}) we get 
\begin{eqnarray*}
\left\| x-\widetilde{x}\right\| _{s,t,\beta } &\leq &k\hat{\rho}_{f}\left\|
y\right\| _{\beta }\Vert x-\tilde{x}\Vert _{s,t,\infty }+\frac{1}{2}\Vert x-%
\tilde{x}\Vert _{s,t,\beta } \\
&&+k\rho _{f}\left\| y-\tilde{y}\right\| _{\beta }+k\rho _{f}^{2}\Vert
y\otimes (y-\tilde{y})\Vert _{\beta }(t-s)^{\beta -2\varepsilon }.
\end{eqnarray*}%
Or 
\begin{eqnarray}
\left\| x-\widetilde{x}\right\| _{s,t,\beta } &\leq &k\hat{\rho}_{f}\left\|
y\right\| _{\beta }\Vert x-\tilde{x}\Vert _{s,t,\infty }  \notag \\
&&\quad +k\rho _{f}\left\| y-\tilde{y}\right\| _{\beta }+k\rho _{f}^{2}\Vert
y\otimes (y-\tilde{y})\Vert _{\beta }(t-s)^{\beta -2\varepsilon }.
\label{3.38}
\end{eqnarray}%
Notice that%
\begin{equation}
\Vert x-\tilde{x}\Vert _{s,t,\infty }\leq \left| x_{s}-\tilde{x}_{s}\right|
+(t-s)^{\beta }\Vert x-\tilde{x}\Vert _{s,t,\beta }.  \label{3.39}
\end{equation}%
Hence, 
\begin{eqnarray*}
\left\| x-\widetilde{x}\right\| _{s,t,\beta } &\leq &k\hat{\rho}_{f}\left\|
y\right\| _{\beta }[\left| x_{s}-\tilde{x}_{s}\right| +(t-s)^{\beta }\Vert x-%
\tilde{x}\Vert _{s,t,\beta }] \\
&&+k\rho _{f}\left\| y-\tilde{y}\right\| _{\beta }+k\rho _{f}^{2}\Vert
y\otimes (y-\tilde{y})\Vert _{\beta }(t-s)^{\beta -2\varepsilon }.
\end{eqnarray*}%
And consequently, 
\begin{equation}
\left\| x-\widetilde{x}\right\| _{s,t,\beta }\leq k\hat{\rho}_{f}\left\|
y\right\| _{\beta }\left| x_{s}-\tilde{x}_{s}\right| +k\rho _{f}\left\| y-%
\tilde{y}\right\| _{\beta }+k\rho _{f}^{2}\Vert y\otimes (y-\tilde{y})\Vert
_{\beta }(t-s)^{\beta -2\varepsilon }.  \label{3.40}
\end{equation}
Substituting (\ref{3.40}) into (\ref{3.39}) yields%
\begin{eqnarray}
\lefteqn{\Vert x-\tilde{x}\Vert _{s,t,\infty } \leq \left| x_{s}-\tilde{x}%
_{s}\right| +(t-s)^{\beta } \Big( k\hat{\rho}_{f}\left\| y\right\| _{\beta }\left| x_{s}-%
\tilde{x}_{s}\right| } \notag \\
&&    +k\rho _{f}\left\| y-\tilde{y}\right\| _{\beta }+k\rho
_{f}^{2}\Vert y\otimes (y-\tilde{y})\Vert _{\beta }(t-s)^{\beta
-2\varepsilon }\Big) .  \label{3.41}
\end{eqnarray}

Suppose that $y=\widetilde{y}$. Then, Equation (\ref{3.41}) implies  that $x=\widetilde{x}$ in a small interval $[0,\delta]$, and by a recursive argument, the uniqueness follows.

Denote $\kappa =  \frac{1}{\beta (y)}\  $ and $%
t_{n}=n\kappa  $. Set 
\begin{equation*}
Z_{n}=\sup_{0\leq s\leq t_{n}}|x_{s}-\tilde{x}_{s}|
\end{equation*}%
Then  inequality \ (\ref{3.41}) states that 
\begin{equation*}
Z_{n+1}\leq (1+k\rho _{f}\kappa  ^{\beta })Z_{n}+k\rho _{f}\kappa
 ^{\beta }\Vert y-\tilde{y}\Vert _{{\beta }}+k\rho _{f}^{2}\kappa
 ^{2{\beta }-2{\varepsilon }}\Vert y\otimes (y-\tilde{y})\Vert _{{%
\beta }}
\end{equation*}%
Therefore 
\begin{eqnarray*}
Z_{T} &\leq &k(1+k\rho _{f}\kappa  ^\beta)^{T/\kappa  }|x_{0}-%
\tilde{x}_{0}|   \\
&&+k\sum_{l=0}^{T/\kappa  }(1+k\rho _{f}\kappa  ^\beta
 )^{l}\left[ \rho _{f}\kappa  ^{\beta }\Vert y-\tilde{y}\Vert _{{%
\beta }}+\rho _{f}^{2}\kappa  ^{2{\beta }-2{\varepsilon }}\Vert
y\otimes (y-\tilde{y})\Vert _{{\beta }}\right].
\end{eqnarray*}
This implies the desired estimate.
\end{proof}

The following corollary is direct consequence of (\ref{3.38}) and (\ref{3.15}%
).

\begin{corollary}
\label{c.4.3} If $f$ is twice continuously differentiable and $f^{\prime
\prime }$ is Lipschitz continuous and if $x$ and $\tilde{x}$ satisfy 
\begin{equation*}
x_{t}=x_{0}+\int_{0}^{t}f(x_{s})dy_{s}\qquad \mathrm{and}\qquad \tilde{x}%
_{t}=\tilde{x}_{0}+\int_{0}^{t}f(\tilde{x}_{s})d\tilde{y}_{s}\,,
\end{equation*}%
then 
\begin{equation}
\Vert x_{t}-\tilde{x}_{t}\Vert _{{\beta }}\leq \ C
\left\{ |x_{0}-\tilde{x}_{0}|+\Vert y-\tilde{y}%
\Vert _{{\beta }}+\Vert y\otimes (y-\tilde{y})\Vert _{{\beta }%
}\right\} \,,  \label{3.30}
\end{equation}%
where we use the notation of Theorem 4.1.
\end{corollary}

 \setcounter{equation}{0}
\section{Stochastic Differential Equations}

Suppose that $B=\{B_{t}=(B_{t}^{1},B_{t}^{2},\ldots ,B_{t}^{d})\}$ is a $d$%
-dimensional Brownian motion. Fix a time interval $[0,T]$. Define%
\begin{equation*}
\left( B\otimes B\right) _{s,t}=\int_{s}^{t}\left( B_{r}-B_{s}\right)
d\circ B_{r},
\end{equation*}%
where the stochastic integral is a Stratonovich integral. That is,

\begin{equation*}
\left( B\otimes B\right) _{s,t}^{i,j}=\left\{ 
\begin{array}{ccc}
\frac{1}{2}(B_{t}^{i}-B_{s}^{i})^{2} & \text{if} & i=j \\ 
\int_{s}^{t}\left( B_{r}^{i}-B_{s}^{i}\right) dB_{r}^j  & \text{if} & i\neq j%
\end{array}%
\right. ,\ 
\end{equation*}%
where the stochastic integral is an It\^{o} integral. It is not difficult to
show that we can choose a version of $\left( B\otimes B\right) _{s,t}$ in 
such a way that $(B,B,B\otimes B)$ constitutes a $\beta $-H\"{o}lder
continuous multiplicative functional, for a fixed $\beta \in (1/3,1/2)$. 

As a first application of Theorem  \ref{t.3.1} and (\ref{eqa}) we deduce that
the Stratonovich stochastic integral $\int_0^T f(B_r) d\circ B_r$
has the following path-wise expression
\begin{eqnarray}
\int_{0}^{T}f(B_{r})\circ dB_{r} &=&(-1)^{\alpha }\ \sum_{i=1}^{d}\int_{0}^{T}%
\widehat{D}_{0+}^{\alpha }f_{i}\left( B\right) _r  (D_{T-}^{1-\alpha
}B_{T-}^{i})_r dr  \notag \\
&&\!\!\!\!  +\sum_{i=1}^{m}\sum_{j=1}^{d}\int_{0}^{T}D_{0+}^{2\alpha -1}\partial
_{i}f_{j}\left( B\right) _r{\Lambda }_{r}^{T}(B^{i}\otimes T^{j})dr\,.
\end{eqnarray}

We can apply Theorem \ref{th4} and deduce the existence   of a
solution for the stochastic differential equation in $\mathbb{R}^{m}$ 
\begin{equation}
X_{t}=X_{0}+\int_{0}^{t}f(X_{s})dB_{s},  \label{4.2}
\end{equation}%
where the initial condition $X_{0}$ is an arbitrary random variable, and the
function $f:\mathbb{R}^{m}\rightarrow \mathbb{R}^{md}$ is a continuously
differentiable function such that $\ f^{\prime }$ is $\lambda $-H\"{o}lder
continuous, where $\lambda >\frac{1}{\beta }-2$, and $f$ and $f^{\prime }$
are bounded. \ By Theorem \ref{th5} the solution is unique if $f$ is twice
continuously differentiable with bounded derivatives and $f^{\prime \prime }$
is $\lambda $-H\"{o}lder continuous, where $\lambda >\frac{1}{\beta }-2.$   \ \  The
stochastic integral here is a path-wise integral which depends on $B$ and $%
B\otimes B$.

We have also the stability type results (\ref{3.15}) and (\ref{3.30}). In
particular, if $B^{{\varepsilon }}$ is a piece-wise smooth approximation of $%
B$ such that 
\begin{equation*}
\Vert B^{{\varepsilon }}-B\Vert _{{\beta }}\qquad \mathrm{and}\qquad \Vert
B\otimes (B^{{\varepsilon }}-B)\Vert _{{\beta }}
\end{equation*}%
converge to zero with a certain rate, then according to Corollary \ref{c.4.3}%
, $\Vert X-X^{{\varepsilon }}\Vert _{{\beta }}$ will also converge to $0$
with the same rate, where 
\begin{equation*}
X_{t}^{{\varepsilon }}=X_{0}+\int_{0}^{t}f(X_{s}^{{\varepsilon }})dB_{s}^{{%
\varepsilon }}.
\end{equation*}%
In particular, this implies that the stochastic process solution of (\ref%
{4.2}) coincides with the solution of the Stratonovich stochastic
differential equation%
\begin{equation}
X_{t}=X_{0}+\int_{0}^{t}f(X_{s})d^{o}B_{s}^{H}.  \label{4.3}
\end{equation}%
In this section we will apply these results in order to obtain the almost
sure rate of convergence \ of the Wong-Zakai approximation to the stochastic
differential equation \ (\ref{4.2}). That is, we will consider the rate of
convergence in H\"{o}lder norm when we approximate the Brownian motion by a
polygonal line.

In order to get a precise rate for these approximations we will make use of
the following exact modulus of continuity of the Brownian motion. There
exists a random variable $G$ such that almost surely for any $s,t\in \lbrack
0,T]$ we have%
\begin{equation}
\left| B_{t}-B_{s}\right| \leq G|t-s|^{1/2}\sqrt{\log \left( |t-s|^{-1}\right) }
\label{4.3b}
\end{equation}

Let $\pi =\left\{ 0=t_{0}<t_{1}<\cdots <t_{n}=T\right\} $ be the uniform
partition of the interval $[0,T]$. That is $t_{k}=\frac{kT}{n}$, $k=0,\ldots
,n$. We denote by $B^{\pi \text{ }}$the polygonal approximation of the
Brownian motion defined by%
\begin{equation*}
B_{t}^{\pi }=\sum_{k=0}^{n-1}\left( B_{t_{k}}+\frac{n}{T}\left(
t-t_{k}\right) \left( B_{t_{k+1}}-B_{t_{k}}\right) \right) \mathbf{1}%
_{(t_{k},t_{k+1}]}(t).
\end{equation*}

We have the following result

\begin{lemma}  \label{lema1}
There exist a random variable $C_{T,\beta}$ such that
\begin{eqnarray}
\Vert B-B^{\pi }\Vert _{{\beta }} &\leq &C_{T,\beta }n^{\beta -1/2}\sqrt{%
\log n}  \label{4.5} \\
\Vert B\otimes (B-B^{\pi })\Vert _{{\beta }} &\leq &C_{T,\beta }n^{\beta
-1/2}\sqrt{\log n}.  \label{4.6}
\end{eqnarray}%
\end{lemma}

\begin{proof}
Fix $0<s<t<T$ and assume that $s\in \lbrack t_{l},t_{l+1}]$ and $t\in
\lbrack t_{k},t_{k+1}]$. Let us first estimate 
\begin{equation*}
h_{1}(s,t)=\frac{1}{(t-s)^{{\beta }}}|B_{t}^{\pi }-B_{t}-(B_{s}^{\pi
}-B_{s})|\,.
\end{equation*}%
If $t-s\geq \frac{T}{n}$, then using (\ref{4.3b}) we obtain 
\begin{eqnarray*}
\left| h_{1}(s,t)\right|  &\leq &T^{-\beta }n^{{\beta }}\left[ \left|
B_{t_{k}}-B_{t}+\frac{n}{T}\left( t-t_{k}\right) \left(
B_{t_{k+1}}-B_{t_{k}}\right) \right| \right.  \\
&&\left. +\left| B_{t_{l}}-B_{s}+\frac{n}{T}\left( s-t_{l}\right) \left(
B_{t_{l+1}}-B_{t_{l}}\right) \right| \right]  \\
&\leq &4GT^{-\beta +1/2}n^{-1/2+\beta }\sqrt{\log \left( n/T\right) }.
\end{eqnarray*}%
If $t-s<\frac{T}{n}$, then there are two cases. Suppose first that $s,t\in
\lbrack t_{k},t_{k+1}]$. In this case, \ if $n$ is large enough ($%
n>Te^{2/(1-2\beta )}$) we obtain using (\ref{4.3b})%
\begin{eqnarray*}
\left| h_{1}(s,t)\right|  &\leq &\frac{|B_{t}-B_{s}|}{(t-s)^{{\beta }}}
+\frac{n}{T} \frac{|B_{t_{k+1}}-B_{t_k}|}{(t-s)^\beta}(t-s)\\
&\leq&
G|t-s|^{\frac{1}{2}-{\beta }}\sqrt{\log |t-s|^{-1}} +GT^{-1/2}\sqrt{\log (n/T)}\ 
n^{1-1/2} (t-s)^{1-\beta} 
\\
&\leq &GT^{-\beta +1/2}n^{-1/2+\beta }\sqrt{\log \left( n/T\right) }.
\end{eqnarray*}%
On the other hand, if $s\in \lbrack t_{k-1},t_{k}]$ and $t\in \lbrack
t_{k},t_{k+1}]$ we have, again if $n$ is large enough 
\begin{eqnarray*}
\left| h_{1}(s,t)\right|  &\leq &\frac{1}{(t-s)^{{\beta }}}\left|
B_{t_{k}}-B_{t}+\frac{n}{T}\left( t-t_{k}\right) \left(
B_{t_{k+1}}-B_{t_{k}}\right) \right.  \\
&&\left. -\left\{ B_{t_{k}}-B_{s}-\frac{n}{T}\left( t_{k}-s\right) \left(
B_{t_{k}}-B_{t_{k-1}}\right) \right\} \right|  \\
&\leq &\frac{1}{(t-s)^{{\beta }}}\left[ |B_{t}-B_{s}|+\frac{n}{T}(t-s)\left(
|B_{t_{k}}-B_{t_{k-1}}|+|B_{t_{k+1}}-B_{t_{k}}|\right) \right]  \\
&\leq &\ \frac{G}{(t-s)^{{\beta }}}\left[ |t-s|^{1/2}\sqrt{\log |t-s|^{-1}}%
+2(t-s)\left( \frac{n}{T}\right) ^{1/2}\sqrt{\log \left( n/T\right) }\right] 
\\
&\leq &3GT^{-\beta +1/2}n^{-1/2+\beta }\sqrt{\log \left( n/T\right) }.
\end{eqnarray*}%
This proves (\ref{4.5}). 

Now we turn to the estimate of the term 
\begin{equation*}
h_{2}(s,t)=\frac{1}{(t-s)^{2{\beta }}}\left|
\int_{s}^{t}(B_{u}^{i}-B_{s}^{i})dB_{u}^{j,\pi
}-\int_{s}^{t}(B_{u}^{i}-B_{s}^{i})dB_{u}^{j}\right| \,
\end{equation*}%
for $i\neq j$ (the case $i=j$ is obvious from (\ref{4.5}). We \ claim that
the there exists a random variable $Z$ such that, almost surely, for all $%
s,t\in \lbrack 0,T]$ we have%
\begin{equation}
\left| \int_{s}^{t}(B_{u}^{i}-B_{s}^{i})dB_{u}^{j}\right| \leq Z|t-s|\log
|t-s|^{-1}.  \label{4.7}
\end{equation}%
In fact, it suffices to show this inequality almost surely for all $s$ and $t
$ rational numbers. If we fix $s$, the process $\left\{ M_{t},t\in \lbrack
s,T]\right\} $%
\begin{equation*}
M_{t}=\int_{s}^{t}(B_{u}^{i}-B_{s}^{i})dB_{u}^{j}
\end{equation*}%
is a continuous martingale and it can be represented as a time-changed
Brownian motion:%
\begin{equation*}
M_{t}=W_{\int_{s}^{t}(B_{u}^{i}-B_{s}^{i})^{2}du}.
\end{equation*}%
As a consequence, applying (\ref{4.3b}) there exists a random variable $G_{1}
$ such that  
\begin{equation*}
|M_{t}|=\left| W_{\int_{s}^{t}(B_{u}^{i}-B_{s}^{i})^{2}du}\right| \leq
G\left( \int_{s}^{t}(B_{u}^{i}-B_{s}^{i})^{2}du\right) ^{1/2}\sqrt{\log
\left( \int_{s}^{t}(B_{u}^{i}-B_{s}^{i})^{2}du\right) ^{-1}}
\end{equation*}%
and  again (\ref{4.3b}), applied to $B_u^i-B_s^i$,   yields%
\begin{equation*}
|M_{t}|\leq G_{1}G_{2}\left( \int_{s}^{t}(u-s)\log |u-s|^{-1}du\right) ^{1/2}%
\sqrt{\log \left( G_{2}^{2}\int_{s}^{t}(u-s)\log |u-s|^{-1}du\right) ^{-1}},
\end{equation*}%
for some random variable $G_{2}$. We have for $|t-s|\leq 1$%
\begin{equation*}
\int_{s}^{t}(u-s)\log |u-s|^{-1}du=|t-s|^{2}\left( \frac{1}{4}+\frac{1}{2}%
\log |t-s|^{-1}\right) ,
\end{equation*}%
and this implies easily the estimate (\ref{4.7}).

Suppose first  that $t-s\geq \frac{T}{n}$. Then%
\begin{eqnarray*}
h_{2}(t,s) &=&\frac{1}{(t-s)^{2\beta }}\left|
\int_{s}^{t}(B_{u}^{i}-B_{s}^{i})dB_{u}^{j,\pi
}-\int_{s}^{t}(B_{u}^{i}-B_{s}^{i})dB_{u}^{j}\right|  \\
&=& \frac{1}{(t-s)^{2\beta }}\left|
(B_t^i-B_s^i) (B_t^{j,\pi} -B_s^{j,\pi} ) - 
\int_{s}^{t}(B_{u}^{j,\pi }-B_{s}^{j,\pi })dB_{u}^i\right.  \\
&& \qquad \left. -   
(B_t^i-B_s^i) (B_t^{j} -B_s^{j} ) +  
\int_{s}^{t}(B_{u}^{j }-B_{s}^{j})dB_{u}^i\right|  \\
 &\leq &\frac{1}{(t-s)^{2\beta }}\left| (B_{t}^{i}-B_{s}^{i})(B_{t}^{j,\pi
}-B_{s}^{j,\pi }-B_{t}^{j}-B_{s}^{j})\right|  \\
&&+\frac{1}{(t-s)^{2\beta }}\left| \int_{s}^{t}\left[ B_{u}^{j,\pi
}-B_{u}^{j}\right] dB_{u}^{i}\right|  \\
&=&A_{1}+A_{2}.
\end{eqnarray*}%
Using (\ref{4.3b}) and (\ref{4.5}) the term $A_{1}$ can be estimated as
follows%
\begin{eqnarray}
A_{1} &\leq &G|t-s|^{1/2-\beta }\sqrt{\log |t-s|^{-1}}\Vert B^{j,\pi
}-B^{j}\Vert _{{\beta }}  \notag \\
&\leq &C_{T,\beta }n^{\beta -1/2}\sqrt{\log n}.  \label{4.10}
\end{eqnarray}%
For the term $A_{2}$ we proceed as in the proof of the estimate (\ref{4.7}).
We have%
\begin{equation*}
\int_{s}^{t}\left[ B_{u}^{j,\pi }-B_{u}^{j}\right] dB_{u}^{i}=W_{%
\int_{s}^{t}\left( B_{u}^{j,\pi }-B_{u}^{j}\right) ^{2}du}
\end{equation*}%
where $W$ is a Brownian motion. \ As a consequence, using that%
\begin{equation*}
\Vert B-B^{\pi }\Vert _{{\infty }}\leq C_{T,\beta }n^{-1/2}\sqrt{\log \left(
n/T\right) }
\end{equation*}
(this estimate is proved as (\ref{4.5})) we get 
\begin{eqnarray}
A_{2} &\leq &\frac{G}{(t-s)^{2\beta }}\left( \int_{s}^{t}(B_{u}^{j,\pi
}-B_{u}^{j})^{2}du\right) ^{1/2}\sqrt{\log \left( \int_{s}^{t}(B_{u}^{j,\pi
}-B_{u}^{j})^{2}du\right) ^{-1}}  \notag \\
&\leq &C_{T,\beta }(t-s)^{1/2-2\beta }n^{-1/2}\sqrt{\log n}\sqrt{\log \left[
(t-s)^{-1}n\left( \log n\right) ^{-1}\right] }  \notag \\
&\leq &C_{T,\beta }n^{\beta -1/2}\sqrt{\log n}.  \label{4.11}
\end{eqnarray}%
Suppose now that   $t-s<\frac{T}{n}$. We make the decomposition%
\begin{eqnarray*}
h_{2}(s,t) &\leq &\frac{1}{(t-s)^{2{\beta }}}\left( \left|
\int_{s}^{t}(B_{u}^{i}-B_{s}^{i})dB_{u}^{j,\pi }\right| +\left|
\int_{s}^{t}(B_{u}^{i}-B_{s}^{i})dB_{u}^{j}\right| \right)  \\
&=&B_{1}+B_{2}.
\end{eqnarray*}%
Then   (\ref{4.7}) yields 
\begin{equation}
B_{2}\leq Z|t-s|^{1-2\beta }\log |t-s|^{-1}\leq C_{T,\beta }n^{2\beta
-1}\log n\,.  \label{4.12}
\end{equation}
In order to handle the term $B_{1}$, assume first that $s,t\in \lbrack
t_{k},t_{k+1}]$. Then 
\begin{equation*}
\int_{s}^{t}(B_{u}^{i}-B_{s}^{i})dB_{u}^{j,\pi }=\frac{n}{T}\left(
B_{t_{k+1}}^{j}-B_{t_{k}}^{j}\right) \int_{s}^{t}(B_{u}^{i}-B_{s}^{i})du
\end{equation*}%
and we obtain%
\begin{equation*}
B_{1}\leq C_{T,\beta }(t-s)^{1-2\beta }\log n\ \leq C_{T,\beta }n^{1-2\beta
}\log n.
\end{equation*}%
Finally, if  $s\in \lbrack t_{k-1},t_{k}]$ and $t\in \lbrack t_{k},t_{k+1}]$
we have%
\begin{eqnarray*}
B_{1} &\leq &(t-s)^{-2\beta }\frac{n}{T}\left| \left(
B_{t_{k}}^{j}-B_{t_{k-1}}^{j}\right)
\int_{s}^{t_{k}}(B_{u}^{i}-B_{s}^{i})du+\left(
B_{t_{k+1}}^{j}-B_{t_{k}}^{j}\right)
\int_{t_{k}}^{t}(B_{u}^{i}-B_{s}^{i})du\right|  \\
&\leq &C_{T,\beta }n^{1-2\beta }\log n.
\end{eqnarray*}
The proof is  now complete.
\end{proof}

As a consequence, we can establish the following result.
\begin{theorem}
Let $f:\mathbb{R}^{m}\rightarrow \mathbb{R}^{md}$ be continuously
differentiable with bounded derivative up to forth order and let $X$ satisfy 
\begin{equation*}
X_{t}=X_{0}+\int_{0}^{t}f(X_{s})dB_{s}
\end{equation*}%
If $X_{t}^{\pi }$ satisfies the following ordinary differential equation 
\begin{equation*}
X_{t}^{\pi }=X_{0}+\int_{0}^{t}f(X_{s}^{\pi })dB_{s}^{\pi }\,,
\end{equation*}%
then for any ${\beta }\in (1/3,1/2)$, there is a random constant $C_{T,\beta
}\in (0,\infty )$ such that 
\begin{equation}
\Vert X-X^{\pi }\Vert _{{\beta }}\leq C_{T,\beta }n^{\beta -1/2}\sqrt{\log n}%
.  \label{4.4}
\end{equation}
\end{theorem}

\begin{proof}
The result is a straightforward consequence of Lemma \ref{lema1}
and Theorem \ref{th5}.  
\end{proof}
 
  \setcounter{equation}{0}
\section{Appendix}

\begin{proof}[Proof of Lemma \ref{l.2.2}]
The fractional integration by parts formula (\ref{1.8}) yields 
\begin{eqnarray*}
&&\int_{a}^{b}d\xi \ \int_{\xi }^{b}\varphi (\xi ,\eta )\frac{\partial
^{2}\psi }{\partial \xi \partial \eta }(\xi ,\eta )d\eta \\
&=&(-1)^{1-\alpha }\int_{a}^{b}d\xi \int_{\xi }^{b}D_{b-}^{\alpha ,\eta
}\varphi _{b-}(\xi ,\cdot )(\eta )D_{\xi +}^{1-\alpha ,\eta }\frac{\partial
\psi }{\partial \xi }(\xi ,\cdot )(\eta )d\eta .
\end{eqnarray*}%
The operators $D_{\xi +}^{1-\alpha ,\eta }$ and $\frac{\partial \psi }{%
\partial \xi }$ commute, as it follows from the following computations:

\begin{eqnarray*}
&&D_{\xi +}^{1-\alpha ,\eta }\frac{\partial \psi }{\partial \xi }(\xi ,\eta )
\\
&=&\frac{1}{\Gamma \left( \alpha \right) }\left( (\eta -\xi )^{\alpha -1}%
\frac{\partial \psi }{\partial \xi }(\xi ,\eta )+(1-\alpha )\int_{\xi
}^{\eta }\frac{\frac{\partial \psi }{\partial \xi }(\xi ,\eta )-\frac{%
\partial \psi }{\partial \xi }(\xi ,\eta ^{\prime })}{\left( \eta -\eta
^{\prime }\right) ^{2-\alpha }}d\eta ^{\prime }\right) \\
&=&\frac{1}{\Gamma (\alpha )}\Bigg\{\frac{\partial }{\partial \xi }\left[
(\eta -\xi )^{\alpha -1}\psi (\xi ,\eta )\right] +(\alpha -1)(\eta -\xi
)^{\alpha -2}\psi (\xi ,\eta ) \\
&&+(1-\alpha )\frac{\partial }{\partial \xi }\int_{\xi }^{\eta }\frac{\psi
(\xi ,\eta )-\psi (\xi ,\eta ^{\prime })}{\left( \eta -\eta ^{\prime
}\right) ^{2-\alpha }}d\eta ^{\prime }+(1-\alpha )\frac{\psi (\xi ,\eta )}{%
\left( \eta -\xi \right) ^{2-\alpha }}\Bigg\} \\
&=&\frac{1}{\Gamma (\alpha )}\frac{\partial }{\partial \xi }\Bigg\{(\eta
-\xi )^{\alpha -1}\psi (\xi ,\eta )+(1-\alpha )\int_{\xi }^{\eta }\frac{\psi
(\xi ,\eta )-\psi (\xi ,\eta ^{\prime })}{\left( \eta -\eta ^{\prime
}\right) ^{2-\alpha }}d\eta ^{\prime }\Bigg\} \\
&=&\frac{\partial }{\partial \xi }D_{\xi +}^{1-\alpha ,\eta }\psi (\xi ,\eta
).
\end{eqnarray*}%
Thus 
\begin{eqnarray*}
&&\int_{a}^{b}d\xi \ \int_{\xi }^{b}\varphi (\xi ,\eta )\frac{\partial
^{2}\psi }{\partial \xi \partial \eta }(\xi ,\eta )d\eta  \\
&=&  (-1)^{1-\alpha}  
 \int_{a}^{b}d\eta
\int_{a}^{\eta }D_{b-}^{\alpha ,\eta }\varphi _{b-}(\xi ,\cdot )(\eta )\frac{%
\partial }{\partial \xi }D_{\xi +}^{1-\alpha ,\eta }\psi (\xi ,\eta )d\xi .
\end{eqnarray*}%
Hence, applying again (\ref{1.8}) we obtain (\ref{1.9}) with ${\Gamma }%
^{\alpha }\psi (\xi ,\eta )$ given by (\ref{1.10}).
\end{proof}

\bigskip \noindent\textbf{Remarks:}

\medskip \noindent \textbf{1. }Formula (\ref{1.9}) holds if $\varphi $ is of
class $C^{2}$ in $a<\xi <\eta <b$ and%
\begin{equation*}
\int_{a}^{b}\int_{a}^{\eta }\left| D_{a+}^{\alpha ,\xi }D_{b-}^{\alpha ,\eta
}\varphi _{b-}(\xi ,\eta )\right| \ d\xi d\eta <\infty .
\end{equation*}

\medskip \noindent\textbf{2. \ }Under the conditions of the above lemma, we
also have ${\Gamma }^{\alpha }\psi (\xi ,\eta )=D_{\xi +}^{1-\alpha ,\eta
}D_{\eta -}^{1-\alpha ,\xi }\psi (\xi ,\eta )$.

\medskip \noindent \textbf{3. \ }The operator ${\Gamma }^{\alpha }$ can also
be expressed as follows.%
\begin{eqnarray*}
&&{\Gamma }^{\alpha }\psi (\xi ,\eta ) \\
&=&\frac{(-1)^{1-\alpha }}{\Gamma \left( \alpha \right) ^{2}}\left\{ (\eta
-\xi )^{\alpha -1}\left( (\eta -\xi )^{\alpha -1}\psi (\xi ,\eta )+(1-\alpha
)\int_{\xi }^{\eta }\frac{\psi (\xi ,\eta )-\psi (\xi ,\eta ^{\prime })}{%
\left( \eta -\eta ^{\prime }\right) ^{2-\alpha }}d\eta ^{\prime }\right)
\right. \\
&&+(1-\alpha )\int_{\xi }^{\eta }\left( \xi ^{\prime }-\xi \right) ^{\alpha
-2}\Bigg[(\eta -\xi )^{\alpha -1}\psi (\xi ,\eta )-(\eta -\xi ^{\prime
})^{\alpha -1}\psi (\xi ^{\prime },\eta ) \\
&&\left. +(1-\alpha )\left( \int_{\xi }^{\eta }\frac{\psi (\xi ,\eta )-\psi
(\xi ,\eta ^{\prime })}{\left( \eta -\eta ^{\prime }\right) ^{2-\alpha }}%
d\eta ^{\prime }-\int_{\xi ^{\prime }}^{\eta }\frac{\psi (\xi ^{\prime
},\eta )-\psi (\xi ^{\prime },\eta ^{\prime })}{\left( \eta -\eta ^{\prime
}\right) ^{2-\alpha }}d\eta ^{\prime }\right) \Bigg]d\xi ^{\prime }\right\}
\\
&=&\frac{(-1)^{1-\alpha }}{\Gamma \left( \alpha \right) ^{2}}\left\{ (\eta
-\xi )^{2\alpha -2}\psi (\xi ,\eta )+(1-\alpha )(\eta -\xi )^{\alpha
-1}\int_{\xi }^{\eta }\frac{\psi (\xi ,\eta )-\psi (\xi ,\eta ^{\prime })}{%
\left( \eta -\eta ^{\prime }\right) ^{2-\alpha }}d\eta ^{\prime }\right. \\
&&+(1-\alpha )\int_{\xi }^{\eta }\frac{(\eta -\xi )^{\alpha -1}\psi (\xi
,\eta )-(\eta -\xi ^{\prime })^{\alpha -1}\psi (\ \xi ^{\prime },\eta )}{%
\left( \xi ^{\prime }-\xi \right) ^{2-\alpha }}d\xi ^{\prime } \\
&&+(1-\alpha )^{2}\int_{\xi }^{\eta }\int_{\xi ^{\prime }}^{\eta }\frac{\psi
(\xi ,\eta )-\psi (\xi ,\eta ^{\prime })-\psi (\xi ^{\prime },\eta )+\psi
(\xi ^{\prime },\eta ^{\prime })}{\left( \xi ^{\prime }-\xi \right)
^{2-\alpha }\left( \eta -\eta ^{\prime }\right) ^{2-\alpha }}d\eta ^{\prime
}d\xi ^{\prime } \\
&&\left. +(1-\alpha )^{2}\int_{\xi }^{\eta }\int_{\xi }^{\xi ^{\prime }}%
\frac{\psi (\xi ,\eta )-\psi (\xi ,\eta ^{\prime })}{\left( \xi ^{\prime
}-\xi \right) ^{2-\alpha }\left( \eta -\eta ^{\prime }\right) ^{2-\alpha }}%
d\eta ^{\prime }d\xi ^{\prime }\right\} \,.
\end{eqnarray*}%
Exchanging the integration order, we see the last double integral equals to 
\begin{equation*}
\int_{\xi }^{\eta }\frac{\psi (\xi ,\eta )-\psi (\xi ,\eta ^{\prime })}{%
\left( \eta -\eta ^{\prime }\right) ^{2-\alpha }}\int_{\eta ^{\prime
}}^{\eta }\frac{1}{\left( \xi ^{\prime }-\xi \right) ^{2-\alpha }}d\xi
^{\prime }d\eta ^{\prime }.
\end{equation*}
This leads to the following expression for $\Gamma ^{\alpha }$
\begin{eqnarray}
\lefteqn{\Gamma ^{\alpha }\psi (\xi ,\eta ) =\frac{(-1)^{1-\alpha }}{\Gamma \left(
\alpha \right) ^{2}}\left\{ (\eta -\xi )^{2\alpha -2}\psi (\xi ,\eta )\right.}
\notag \\
&&\quad +(1-\alpha )\int_{\xi }^{\eta }\frac{(\eta -\xi )^{\alpha -1}\psi
(\xi ,\eta )-(\eta -\xi ^{\prime })^{\alpha -1}\psi (\ \xi ^{\prime },\eta )%
}{\left( \xi ^{\prime }-\xi \right) ^{2-\alpha }}d\xi ^{\prime }  \notag \\
&&\quad +(1-\alpha )^{2}\int_{\xi }^{\eta }\int_{\xi ^{\prime }}^{\eta }%
\frac{\psi (\xi ,\eta )-\psi (\xi ,\eta ^{\prime })-\psi (\xi ^{\prime
},\eta )+\psi (\xi ^{\prime },\eta ^{\prime })}{\left( \xi ^{\prime }-\xi
\right) ^{2-\alpha }\left( \eta -\eta ^{\prime }\right) ^{2-\alpha }}d\eta
^{\prime }d\xi ^{\prime }  \notag \\
&&\quad +\left. (\alpha -1)\int_{\xi }^{\eta }\frac{\phi (\xi ,\eta )-\phi
(\xi ,\eta ^{\prime })}{(\eta -\eta ^{\prime })^{2-\alpha }}(\eta ^{\prime
}-\xi )^{\alpha -1}d\eta ^{\prime }\right\} \,.  \label{1.11}
\end{eqnarray}

Consider the kernel $K_{s,t}(\xi ,\eta )$ defined in (\ref{2.7}), that is,%
\begin{equation*}
K_{s,t}(\xi ,\eta )=D_{s+}^{1,{\alpha -{\varepsilon }}}D_{t-}^{2,{\alpha -{%
\varepsilon }}}G_{t-}(s,\xi ,\eta ),
\end{equation*}%
where%
\begin{eqnarray}
G(s,\xi ,\eta ) &=&C_{\alpha }\left( \xi -s\right) ^{\alpha -1}\left( \eta
-\xi \right) ^{\alpha -1}\int_{0}^{1}q^{2\alpha -2}(1-q)^{-\alpha }(1+(1-q)%
\frac{\xi -s}{\eta -\xi })^{-1}dq  \notag \\
&=&\left( \xi -s\right) ^{\alpha -1}\left( \eta -\xi \right) ^{\alpha
-1}\phi \left( \frac{\xi -s}{\eta -\xi }\right) ,  \label{e.7.1}
\end{eqnarray}%
\begin{equation}
\phi (z)=C_{\alpha }\int_{0}^{1}q^{2\alpha -2}(1-q)^{-\alpha
}(1+(1-q)z)^{-1}dq=C_{\alpha }\int_{0}^{1}(1-q)^{2\alpha -2}q^{-\alpha
}(1+qz)^{-1}dq,  \label{e.7.2}
\end{equation}%
and $C_{\alpha }$ is given as the coefficient in (\ref{2.6}).

\begin{lemma}
\label{l.7.1} Let $1/2<\alpha<1$. The function $\phi (z)$ defined in (\ref%
{e.7.2}) satisfies $\phi (0)<\infty $, $\phi $ is decreases to zero as $z$
tends to infinity. If $\beta <1-\alpha $, then 
\begin{equation}
\phi (z)\leq cz^{-\beta }  \label{e.7.3}
\end{equation}
Moreover, if $\beta <2-\alpha $, 
\begin{equation}
\left| \phi ^{\prime }(z)\right| \leq cz^{-\beta }.  \label{e.7.4}
\end{equation}
\end{lemma}

\begin{lemma}
\label{l.7.2} \bigskip The kernel $K_{s,t}(\xi ,\eta )$ satisfies%
\begin{equation}
\sup_{0\le s<t\le T }\int_{s< \xi < \eta < t}\left| K_{s,t}(\xi ,\eta )\right|
d\xi d\eta <\infty .  \label{e.7.5}
\end{equation}
\end{lemma}

\begin{proof}
To simplify the notation we omit the dependence on the variable $s$ in $%
G(s,\xi ,\eta )$. Also, $c$ will denote a generic constant depending on $%
\alpha $ and $\varepsilon $. We have 
\begin{eqnarray*}
&&D_{s+}^{1,{\alpha -{\varepsilon}} }D_{t-}^{2,{\alpha -{\varepsilon}}
}G_{t-}(\xi ,\eta ) \\
&=&cD_{s+}^{1,{\alpha -{\varepsilon}} }\left( \frac{G(\xi ,\eta )-G(\xi ,t)}{%
(t-\eta )^{{\alpha -{\varepsilon}} }}+({\alpha -{\varepsilon}}) \int_{\eta
}^{t}\frac{G(\xi ,\eta )-G(\xi ,\eta ^{\prime })}{(\eta ^{\prime }-\eta )^{{%
\alpha -{\varepsilon}} +1}}d\eta ^{\prime }\right) \\
&=&c\left( \frac{G(\xi ,\eta )-G(\xi ,t)}{(t-\eta )^{{\alpha -{\varepsilon}}
}(\xi -s)^{{\alpha -{\varepsilon}} }}+\frac{1}{(\xi -s)^{{\alpha -{%
\varepsilon}} }}\int_{\eta }^{t}\frac{G(\xi ,\eta )-G(\xi ,\eta ^{\prime })}{%
(\eta ^{\prime }-\eta )^{{\alpha -{\varepsilon}} +1}}d\eta ^{\prime }\right.
\\
&&+({\alpha -{\varepsilon}}) \int_s^{\xi }\frac{G(\xi ,\eta )-G(\xi
,t)-G(\xi ^{\prime },\eta )+G(\xi ^{\prime },t)}{(t-\eta )^{{\alpha -{%
\varepsilon}} }(\xi -\xi ^{\prime })^{{\alpha -{\varepsilon}} +1}}d\xi
^{\prime } \\
&&\left. +({\alpha -{\varepsilon}} ) ^{2}\int_s^{\xi }\int_{\eta }^{t}\frac{%
G(\xi ,\eta )-G(\xi ,\eta ^{\prime })-G(\xi ^{\prime },\eta )+G(\xi ^{\prime
},\eta ^{\prime })}{(\eta ^{\prime }-\eta )^{{\alpha -{\varepsilon}} +1}(\xi
-\xi ^{\prime })^{{\alpha -{\varepsilon}} +1}}d\eta ^{\prime }d\xi ^{\prime
}\right) .
\end{eqnarray*}%
Set%
\begin{eqnarray*}
A_{1} &=&\frac{G(\xi ,\eta )-G(\xi ,t)}{(t-\eta )^{{\alpha -{\varepsilon}}
}(\xi -s)^{{\alpha -{\varepsilon}} }} \\
A_{2} &=&\frac{1}{(\xi -s)^{{\alpha -{\varepsilon}} }}\int_{\eta }^{t}\frac{%
G(\xi ,\eta )-G(\xi ,\eta ^{\prime })}{(\eta ^{\prime }-\eta )^{{\alpha -{%
\varepsilon}} +1}}d\eta ^{\prime } \\
A_{3} &=&\int_s^{\xi }\frac{G(\xi ,\eta )-G(\xi ,t)-G(\xi ^{\prime },\eta
)+G(\xi ^{\prime },t)}{(t-\eta )^{{\alpha -{\varepsilon}} }(\xi -\xi
^{\prime })^{{\alpha -{\varepsilon}} +1}}d\xi ^{\prime } \\
A_{4} &=&\int_s^{\xi }\int_{\eta }^{t}\frac{G(\xi ,\eta )-G(\xi ,\eta
^{\prime })-G(\xi ^{\prime },\eta )+G(\xi ^{\prime },\eta ^{\prime })}{(\eta
^{\prime }-\eta )^{{\alpha -{\varepsilon}} +1}(\xi -\xi ^{\prime })^{{\alpha
-{\varepsilon}} +1}}d\eta ^{\prime }d\xi ^{\prime }.
\end{eqnarray*}%
It suffices to show that%
\begin{equation}
\sup_{0\le s<t\le T}\int_{s< \xi < \eta < t}|A_{i}| d\xi d\eta <\infty \ 
\label{e.7.6}
\end{equation}%
for $i=1,2,3,4.$

\textbf{Step 1} Suppose $i=1$. Using the fact that the function $\phi $ is
bounded we obtain 
\begin{equation*}
G(\xi ,\eta )\leq c\left( \xi -s\right) ^{\alpha -1}\left( \eta -\xi \right)
^{\alpha -1}.
\end{equation*}%
Hence,%
\begin{eqnarray*}
\left| A_{1}\right| &=&\frac{|G(\xi ,\eta )|+ |G(\xi ,t)|}{(t-\eta )^{{%
\alpha -{\varepsilon}} }(\xi -s)^{{\alpha -{\varepsilon}} }} \\
&\leq & c\left( \xi -s\right) ^{\varepsilon -1}\left( \eta -\xi \right)
^{\alpha -1}(t-\eta )^{-\alpha +\varepsilon } \\
&&\quad + c(\xi-s)^{{\varepsilon}-1}(t-\eta)^{{\varepsilon}-1}\,.
\end{eqnarray*}
and (\ref{e.7.6} ) holds for $i=1$.

\textbf{Step 2} Suppose $i=2$. \ We have%
\begin{equation*}
G(\xi ,\eta )-G(\xi ,\eta ^{\prime })=\int_{\eta }^{\eta ^{\prime }}\frac{%
\partial G}{\partial y}(\xi ,y)dy,\ 
\end{equation*}%
and 
\begin{eqnarray*}
\frac{\partial G}{\partial \eta }(\xi ,\eta ) &=&(\alpha -1) (\xi
-s)^{\alpha -1}(\eta -\xi )^{\alpha -2}\phi \left( \frac{\xi -s}{\eta -\xi }%
\right) \\
&&- (\xi -s)^{\alpha -1}(\eta -\xi )^{\alpha -1}\phi ^{\prime }\left( \frac{
\xi -s}{\eta -\xi }\right) \frac{\xi -s}{\left( \eta -\xi \right) ^{2}} \\
&=&(\xi -s)^{\alpha -1}(\eta -\xi )^{\alpha -2}\chi \left( \frac{\xi -s}{%
\eta -\xi }\right) ,
\end{eqnarray*}%
where%
\begin{equation*}
\chi (z)=(\alpha -1) \phi (z)- z\phi ^{\prime }(z).
\end{equation*}%
Notice that, by Lemma \ref{l.7.1} the function $\chi (z)$ is uniformly
bounded. Hence,%
\begin{eqnarray*}
|A_{2}| &\leq &c(\xi -s)^{\varepsilon -1}\int_{\eta }^{t}\int_{\eta }^{\eta
^{\prime }}(y-\xi )^{\alpha -2}(\eta ^{\prime }-\eta )^{-\alpha +\varepsilon
-1}dyd\eta ^{\prime } \\
&\leq &c(\xi -s)^{\varepsilon -1}(\eta -\xi )^{\frac{\varepsilon }{4}%
-1}\int_{\eta }^{t} \int_{\eta }^{\eta ^{\prime }}(\eta ^{\prime }-y)^{\frac{%
\varepsilon }{2}-1}(y-\eta )^{\frac{\varepsilon }{4}-1}dyd\eta ^{\prime } \\
&\leq &c(\xi -s)^{\varepsilon -1}(\eta -\xi )^{\frac{\varepsilon }{4}%
-1}\int_{\eta }^{t} (\eta ^{\prime }-\eta )^{\frac{3\varepsilon }{4}-1}d\eta
^{\prime },
\end{eqnarray*}%
which implies that (\ref{e.7.6} ) holds for $i=2$.$\ $

\textbf{Step 3} Suppose $i=3$. We have%
\begin{equation*}
G(\xi ,\eta )-G(\xi ^{\prime },\eta )=\int_{\xi ^{\prime }}^{\xi }\frac{%
\partial G}{\partial x}(x,\eta )dx,
\end{equation*}%
and 
\begin{eqnarray*}
\frac{\partial G}{\partial \xi }(\xi ,\eta ) &=&(\alpha -1) (\xi -s)^{\alpha
-2}(\eta -\xi )^{\alpha -1}\phi \left( \frac{\xi -s}{\eta -\xi }\right) \\
&&-(\alpha -1) (\xi -s)^{\alpha -1}(\eta -\xi )^{\alpha -2}\phi \left( \frac{%
\xi -s}{\eta -\xi }\right) \\
&&+ (\xi -s)^{\alpha -1}(\eta -\xi )^{\alpha -1}\phi ^{\prime }\left( \frac{%
\xi -s}{\eta -\xi }\right) \frac{\eta -s}{\left( \eta -\xi \right) ^{2}} \\
&=&(\xi -s)^{\alpha -2}(\eta -\xi )^{\alpha -1}(\alpha -1) \phi \left( \frac{%
\xi -s}{\eta -\xi }\right) \\
&&+(\xi -s)^{\alpha -1}(\eta -\xi )^{\alpha -2}\gamma \left( \frac{\xi -s}{%
\eta -\xi }\right) ,
\end{eqnarray*}%
where%
\begin{equation*}
\gamma (z)=(1-\alpha ) \phi (z)+\left( 1+z\right) \phi ^{\prime }(z).
\end{equation*}%
By Lemma \ref{l.7.1} the function $\gamma $ is uniformly bounded. Hence,%
\begin{eqnarray*}
|A_{3}| &\leq &c(t-\eta )^{-\alpha +\varepsilon }\int_s^{\xi }\int_{\xi
^{\prime }}^{\xi }(\xi -\xi ^{\prime })^{-\alpha +\varepsilon -1} \\
&&\times \left[ (x-s)^{\alpha -2}(\eta -x)^{\alpha -1}+(x-s)^{\alpha
-1}(\eta -x)^{\alpha -2}\right] dxd\xi ^{\prime } \\
&\leq &c(t-\eta )^{-\alpha +\varepsilon }(\eta -\xi )^{\alpha -1}\int_s^{\xi
}\int_{\xi ^{\prime }}^{\xi }(\xi -x)^{\frac{\varepsilon }{2}-1}(x-\xi
^{\prime })^{\frac{\varepsilon }{4}-1}(\xi ^{\prime }-s)^{\frac{\varepsilon 
}{4}-1}dxd\xi ^{\prime } \\
&&+c(t-\eta )^{-\alpha +\varepsilon }(\eta -\xi )^{\frac{\varepsilon }{4}%
-1}\int_s^{\xi }\int_{\xi ^{\prime }}^{\xi }(\xi -x)^{\frac{\varepsilon }{4}%
-1}(x-\xi ^{\prime })^{\frac{\varepsilon }{2}-1}(\xi ^{\prime }-s)^{\alpha
-1}dxd\xi ^{\prime } \\
&= &c(t-\eta )^{-\alpha +\varepsilon }(\eta -\xi )^{\alpha -1}(\xi-s)^{{%
\varepsilon}-1}+c (t-\eta )^{-\alpha +\varepsilon }(\eta -\xi )^{\frac{%
\varepsilon }{4} -1} (\xi-s)^{\frac34 {\varepsilon}+\alpha-1}\,.
\end{eqnarray*}%
which implies that (\ref{e.7.6} ) holds for $i=3$.

\textbf{Step 4} Suppose $i=4$. We are going to use the following
decomposition%
\begin{eqnarray*}
&&G(\xi ,\eta )-G(\xi ,\eta ^{\prime })-G(\xi ^{\prime },\eta )+G(\xi
^{\prime },\eta ^{\prime }) \\
&=&-\int_{\eta }^{\eta ^{\prime }}\int_{\xi ^{\prime }}^{\xi }\frac{\partial
^{2}G}{\partial x\partial y}(x,y)dxdy
\end{eqnarray*}%
We need to compute the second derivative:%
\begin{eqnarray*}
\frac{\partial ^{2}G}{\partial \xi \partial \eta } &=&(\xi -s)^{\alpha
-2}(\eta -\xi )^{\alpha -2}\left( (\alpha -1)^{2}-(\alpha -1)(\alpha -2)%
\frac{\xi -s}{\eta -\xi }\right) \phi \left( \frac{\xi -s}{\eta -\xi }\right)
\\
&&-(\xi -s)^{\alpha -1}(\eta -\xi )^{\alpha -1}\left[ \phi ^{\prime \prime
}\left( \frac{\xi -s}{\eta -\xi }\right) \frac{(\eta -s)(\xi -s)}{\left(
\eta -\xi \right) ^{4}}\right. \\
&&\left. +\phi ^{\prime }\left( \frac{\xi -s}{\eta -\xi }\right) \frac{(\eta
-\xi )+2(\xi -s)}{\left( \eta -\xi \right) ^{3}}\right] .
\end{eqnarray*}%
Hence, we can write%
\begin{equation*}
\frac{\partial ^{2}G}{\partial \xi \partial \eta }=(\xi -s)^{\alpha -2}(\eta
-\xi )^{\alpha -2}\psi \left( \frac{\xi -s}{\eta -\xi }\right) ,
\end{equation*}%
where%
\begin{eqnarray*}
\psi (z) &=&\left( (\alpha -1)^{2}-(\alpha -1)(\alpha -2)z\right) \phi
\left( z\right) \\
&&-\phi ^{\prime \prime }\left( z\right) z^{2}(1+z)-\phi ^{\prime
}(z)z(1+2z).
\end{eqnarray*}%
We are going to use the decomposition%
\begin{equation*}
\psi (z)=\psi _{1}(z)+\psi _{2}(z),
\end{equation*}%
where%
\begin{eqnarray*}
\psi _{1}(z) &=&-(\alpha -1)(\alpha -2)z\phi \left( z\right) \\
\psi _{2}(z) &=&(\alpha -1)^{2}\phi \left( z\right) -\phi ^{\prime \prime
}\left( z\right) z^{2}(1+z)-\phi ^{\prime }(z)z(1+2z).
\end{eqnarray*}%
This leads to 
\begin{equation*}
\int_{s\leq \xi \leq \eta \leq t}A_{4}d\xi d\eta =B_{1}+B_{2},
\end{equation*}%
where%
\begin{eqnarray*}
B_{i} &=&-\int_{D}(\eta ^{\prime }-\eta )^{-\alpha +\varepsilon -1}(\xi -\xi
^{\prime })^{-\alpha +\varepsilon -1} \\
&&\times (x-s)^{\alpha -2}(y-\xi )^{\alpha -2}\psi _{i}(\frac{x-s}{y-x})d\xi
^{\prime }dxd\xi d\eta dyd\eta ^{\prime },
\end{eqnarray*}%
$i=1,2$, and%
\begin{equation*}
D=\{(\xi ^{\prime },x,\xi ,\eta ,y,\eta ^{\prime }):s<\xi ^{\prime }<x<\xi
<\eta <y<\eta ^{\prime }<t\}.
\end{equation*}
\textbf{Step 5} Estimation of $B_{1}$. Denote 
\begin{equation*}
D_{1}=\{(\xi ^{\prime },x,\xi ,\eta ,\eta ^{\prime }):s<\xi ^{\prime }<x<\xi
<\eta <\eta ^{\prime }<t\}.
\end{equation*}
Using (\ref{e.7.3}) with $\beta =1-\alpha -\delta $ with $\delta
<\varepsilon/3 $, we obtain%
\begin{eqnarray*}
|B_{1}| &\leq &c\int_{D}(\eta ^{\prime }-\eta )^{-\alpha +\varepsilon
-1}(\xi -\xi ^{\prime })^{-\alpha +\varepsilon -1}(x-s)^{\alpha
-1}(y-x)^{\alpha -3} \\
&&\times \phi \left(\frac{x-s}{y-x} \right) d\xi ^{\prime }dxd\xi d\eta
dyd\eta ^{\prime } \\
&\leq &c\int_{D}(\eta ^{\prime }-\eta )^{-\alpha +\varepsilon -1}(\xi -\xi
^{\prime })^{-\alpha +\varepsilon -1}(x-s)^{2\alpha -2+\delta
}(y-x)^{-2-\delta }d\xi ^{\prime }dxd\xi d\eta dyd\eta ^{\prime } \\
&=&c\int_{D_{1}}(\eta ^{\prime }-\eta )^{-\alpha +\varepsilon -1}(\xi -\xi
^{\prime })^{-\alpha +\varepsilon -1}(x-s)^{2\alpha -2+\delta } \\
&&\times \left[ (\eta -x)^{-1-\delta }-(\eta ^{\prime }-x)^{-1-\delta }%
\right] d\xi ^{\prime }dxd\xi d\eta d\eta ^{\prime } \\
&\leq & c\int_{D_1}(\eta ^{\prime }-\eta )^{-\alpha +\varepsilon -1 }(\xi
-\xi ^{\prime })^{-\alpha +\varepsilon -1}(x-s)^{2\alpha -2+\delta } (\eta
-x)^{-1-\delta } (\eta-\eta^{\prime})^\alpha d\xi ^{\prime }dxd\xi d\eta
d\eta ^{\prime } \\
&\leq &c\int_{D_1}(\eta ^{\prime }-\eta )^{ \varepsilon -1 } (\xi
-x)^{-1+3\delta+\alpha } (x -\xi ^{\prime })^{-2\alpha+{\varepsilon}-3\delta
} \\
&&\quad (x-\xi^{\prime})^{2\alpha-1} (\xi^{\prime}-s)^{-1+\delta }
(\eta-\xi)^{-1+\delta} (\xi -x)^{-2\delta -\alpha }d\xi ^{\prime }dxd\xi
d\eta d\eta ^{\prime } \\
&\leq &c\int_{D_1}(\eta ^{\prime }-\eta )^{\varepsilon -1 } (\xi
-x)^{-1+\delta} (x -\xi ^{\prime })^{{\varepsilon}-3\delta-1 }
(\xi^{\prime}-s)^{-1+\delta } (\eta-\xi)^{-1+\delta} d\xi ^{\prime }dxd\xi
d\eta d\eta ^{\prime } \\
&<&\infty .
\end{eqnarray*}

\textbf{Step 6} Estimation of $B_{2}$. Let us compute the function $\psi
_{2}(z)$:%
\begin{eqnarray*}
\psi _{2}(z) &=&\int_{0}^{1}(1-q)^{2\alpha -2}q^{-\alpha }\left[ (\alpha
-1)^{2}(1+qz)^{-1}\right. \\
&&\left. -2q^{2}(1+qz)^{-3}z^{2}(1+z)+q(1+qz)^{-2}z(1+2z)\right] dq \\
&=&\int_{0}^{1}(1-q)^{2\alpha -2}q^{-\alpha }(1+qz)^{-3}\left[ (\alpha
-1)^{2}(1+qz)^{2}\right. \\
&&\qquad +qz(1+(2-q)z)]dq
\end{eqnarray*}%
This implies that the function $\psi _{2}(z)$ is uniformly bounded. As a
consequence, we deduce the following estimates 
\begin{eqnarray*}
|B_{2}| &\leq &c\int_{D}(\eta ^{\prime }-\eta )^{-\alpha +\varepsilon
-1}(\xi -\xi ^{\prime })^{-\alpha +\varepsilon -1}(x-s)^{\alpha -2}(y-\xi
)^{\alpha -2}d\xi ^{\prime }dxd\xi d\eta dyd\eta ^{\prime } \\
&\leq &c\int_{D}(\eta ^{\prime }-y)^{\frac{\varepsilon }{2}-1}(y-\eta
)^{-\alpha +\frac{\varepsilon }{2}}(\xi -\xi ^{\prime })^{-\alpha
+\varepsilon -1}(x-s)^{\alpha -2}(y-\xi )^{\alpha -2}d\xi ^{\prime }dxd\xi
d\eta dyd\eta ^{\prime } \\
&\leq &c\int_{D}(\eta ^{\prime }-y)^{\frac{\varepsilon }{2}-1}(y-\eta )^{-1+%
\frac{\varepsilon }{4}}(\xi -x)^{\frac{\varepsilon }{2}-1}(x-\xi ^{\prime
})^{-\alpha +\frac{\varepsilon }{2}} \\
&&\times (x-s)^{\alpha -2}(\eta -\xi )^{-1+\frac{\varepsilon }{4}}d\xi
^{\prime }dxd\xi d\eta dyd\eta ^{\prime } \\
&\leq &c\int_{D}(\eta ^{\prime }-y)^{\frac{\varepsilon }{2}-1}(y-\eta )^{-1+%
\frac{\varepsilon }{4}}(\xi -x)^{\frac{\varepsilon }{2}-1}(x-\xi ^{\prime
})^{-1+\frac{\varepsilon }{4}} \\
&&\times (\xi ^{\prime }-s)^{-1+\frac{\varepsilon }{4}}(\eta -\xi )^{-1+%
\frac{\varepsilon }{4}}d\xi ^{\prime }dxd\xi d\eta dyd\eta ^{\prime } \\
&<&\infty .
\end{eqnarray*}
\end{proof}

\end{document}